\documentclass[11pt,english]{article}
\usepackage[T1]{fontenc}
\usepackage[utf8]{inputenc}
\usepackage[letterpaper]{geometry}
\usepackage{amsmath}
\usepackage{graphicx}
\usepackage{amssymb}

\makeatletter
\newtheorem{theorem}{Theorem}[section]

\newtheorem{lemma}[theorem]{Lemma}
\newtheorem{proposition}[theorem]{Proposition}

\numberwithin{equation}{section}

\makeatother

\newcommand{\eqnsection}{
\renewcommand{\theequation}{\thesection.\arabic{equation}}
 \makeatletter   \csname  @addtoreset\endcsname{equation}{section}
   \makeatother}

\def\proof{\noindent\textbf{Proof: $\,$} }
\def\qed{$\Box $}

\def\R{\mathbb{R}}
\def\N{\mathbb{N}}

\def\E{\mathbb{E}}
\def\P{\mathbb{P}}
\def\0{\mathbf{0}}
\def\1{\mathbf{1}}

\def\Var{{\mathop {{\rm Var\, }}}}
\def\Cov{{\mathop {{\rm Cov\, }}}}

\makeatother

\usepackage{babel}

\makeatother

\begin{document}

\title{Large deviations for local times and intersection local times of
fractional Brownian motions and Riemann-Liouville processes}

\author{Xia Chen%
\thanks{Supported in part by NSF grant DMS-0704024.%
} \, \, Wenbo V. Li%
\thanks{Supported in part by NSF grant DMS--0805929, NSFC-6398100, CAS-2008DP173182.%
} \, \, Jan Rosi\'nski%
\thanks{Supported in part by NSA grant MSPF-50G-049.%
}\, \, Qi-Man Shao%
\thanks{Supported in part by Hong Kong RGC CERG 602206 and 602608.%
}}

\date{May, 2010}

\maketitle

\begin{abstract}
In this paper we prove exact forms of large deviations for local times and intersection local times of fractional Brownian motions and Riemann--Liouville processes. We also show that a fractional Brownian motion and the related Riemann--Liouville process behave like constant multiples of each other with regard to large deviations for their local and intersection local times. As a consequence of our large deviation estimates, we derive laws of iterated logarithm for the corresponding local times. The key points of our methods: (1) logarithmic {\em superadditivity} of a normalized sequence of moments of exponentially randomized local time of a fractional Brownian motion; (2) logarithmic {\em subadditivity} of a normalized sequence of moments of exponentially randomized intersection local time   of Riemann--Liouville processes;  (3) comparison of local and intersection  local times based on embedding of a part of a fractional Brownian motion into the reproducing kernel Hilbert space of the Riemann--Liouville process.
\end{abstract}
\bigskip

\begin{quote} {\footnotesize

}
\end{quote}

%60G22 Fractional processes, including fractional Brownian motion;  60J55  Local time and additive functionals;  60F10  Large deviations; 60G15  Gaussian processes;  60G18 Self-similar processes;

\newpage

\baselineskip4pt
\parskip4pt

%\tableofcontents

\baselineskip16pt
\parskip6pt

\section{Introduction } \label{s:intro}

Let $B^{H}(t)$, $t\ge0$ be a standard $d$-dimensional fractional
Brownian motion with index $H\in(0,1)$. That is, $B^{H}(t)$ is a
zero-mean Gaussian process with stationary increments and covariance
function
\[
\E\left[B^{H}(t)B^{H}(s)^{\top}\right]=\frac{1}{2}\left\{ |t|^{2H}+|s|^{2H}-|t-s|^{2H}\right\} I_{d},
\]
 where $I_{d}$ is the identity matrix of size $d$. $B^{H}(t)$ is also
a self-similar process with index $H$. The local time $L_{t}^{x}(B^{H})$
of $B^{H}(t)$ at $x\in\R^{d}$ is defined heuristically as
\[
L_{t}^{x}(B^{H})=\int_{0}^{t}\delta_{x}\left(B^{H}(s)\right) \, ds,\quad t\ge0.
\]
 It is known that $L_{t}^{x}(B^{H})$ exists and is jointly continuous
in $(t,x)$ as long as $Hd<1$. By the self-similarity of a fractional
Brownian motion,
$
L_{t}^{x}(B^{H})\,\overset{d}{=}t^{1-Hd}L_{1}^{x/t^{H}}(B^{H}).
$
In particular,
\begin{align}
L_{t}^{0}(B^{H})\,\overset{d}{=}\, t^{1-Hd}L_{1}^{0}(B^{H}) & .\label{eq:LB-ss}
\end{align}

Our first goal is to investigate large deviations associated with
tail probabilities of $L_{t}^{0}(B^{H})$. By the scaling given
above, we may consider only $t=1$. In the classical case, when $H=1/2$
and $d=1$, it is well known, see the book of Revuz and Yor \cite{RY}, p240, 
that $L_{1}^{0}(B^{1/2})\overset{d}{=}|U|$
with $U\sim\mathcal{N}(0,1)$. Consequently,
\[
\lim_{a\to\infty}a^{-2}\log\mathbb{P}\left\{ L_{1}^{0}(B^{1/2})\ge a\right\} =-\frac{1}{2}.
\]
 In Theorem \ref{th:1} we prove that for a fractional Brownian motion a nontrivial limit
\[
\lim_{a\to\infty}a^{-1/Hd}\log\P\{L_{1}^{0}(B^{H})\ge a\}\]
 exists and we give bounds for this limit. \medskip{}

Closely related to the fractional Brownian motion is the Riemann–Liouville
process $W^{H}(t)$ with index $H>0$ which is defined as a stochastic convolution

\begin{equation}
W^{H}(t)=\int_{0}^{t}(t-s)^{H-1/2}\, dB(s)\,,\quad t\ge0,\label{eq:W}
\end{equation}
 where $B(t)$ is a $d$-dimensional standard Brownian
motion. $\{W^{H}(t)\}_{t\ge0}$
is a self-similar zero-mean Gaussian process with index $H$, as is
$B^{H}(t)$, but $W^{H}(t)$ does not have stationary increments
and there is no upper bound restriction on index $H>0$.
If $L_{t}^{0}(W^{H})$ denotes the local time of $W^{H}(t)$ at 0,
then by the self-similarity we also have
\begin{equation}
L_{t}^{0}(W^{H})\,\overset{d}{=}\, t^{1-Hd}L_{1}^{0}(W^{H})\,.
\label{eq:LW-ss}
\end{equation}
 The relation between $W^{H}(t)$ and $B^{H}(t)$ becomes transparent
when we write a moving average representation of $B^{H}(t)$, $t\in\R,$
in the form
\begin{equation}
B^{H}(t)=c_{H}\int_{-\infty}^{t}\left[(t-s)^{H-1/2}-(-s)_{+}^{H-1/2}\right]\, dB(s),\label{eq:B}
\end{equation}
where
\begin{equation}
c_{H}=\sqrt{2H}\,2^{H}B\left(1-H,H+1/2\right)^{-1/2}\,,\label{eq:c}
\end{equation}
and $B(\cdot,\cdot)$ denotes the beta function. The analytic derivation of $c_{H}$  is given for completeness in the Appendix (a different but equivalent form of $c_{H}$ is also derived in Mishura, \cite[Lemma A.0.1]{M}, by a Fourier analytic method). From \eqref{eq:B} we have a decomposition
\begin{equation}
c_{H}^{-1}B^{H}(t)=W^{H}(t)+Z^{H}(t),\label{eq:dec}
\end{equation}
 where
 \begin{equation}
Z^{H}(t)=\int_{-\infty}^{0}\left[(t-s)^{H-1/2}-(-s)^{H-1/2}\right]\, dB(s)\label{eq:rem}
\end{equation}
 is a process independent of $W^{H}(t)$.

This moving average representation for fractional Brownian motion
was introduced in the pioneering work of Mandelbrot and Van Ness \cite{MV} and used extensively by many authors, sometimes with different normalizing constant $c_H$ in (\ref{eq:c}) (e.g., Li and Linde \cite{LL98} uses $\Gamma(H+1/2)^{-1}$ for $c_H$).

We will show that paths
of $Z^{H}(t),$ away from $t=0,$ can be matched with functions in
the reproducing kernel Hilbert space of $W^{H}(t)$ (Proposition \ref{prop:Z_a},
Section \ref{sub:remainder}). This and the independence of $Z^{H}(t)$
from $W^{H}(t)$ will allow us to show that large deviation constants
of tail probabilities of $L_{1}^{0}(W^{H})$ and $ $ of $L_{1}^{0}(c_{H}^{-1}B^{H})=c_{H}^{d}L_{1}^{0}(B^{H})$
are the same (Theorem \ref{th:2}). In this context we also want to
mention Theorem 3.22 of Xiao, \cite{Xiao-1}, who established bounds
for tail probabilities of the local time $L_{1}^{0}$ of the general
Gaussian processes in the form
\[
-C_{1}\le\liminf_{a\to\infty}\frac{1}{\phi(a)}\log\{L_{1}^{0}\ge a\}\le\limsup_{a\to\infty}\frac{1}{\phi(a)}\log\{L_{1}^{0}\ge a\}\le-C_{2}
\]
 and raised a question on the existence of the limit (Question 3.25,
\cite{Xiao-1}). Further, we cite the paper by Baraka, Mountford and
Xiao (\cite{BMX}) for some similar tail estimate of the local time of
multi-parameter fractional Brownian motions.

Next we will consider $p$ independent copies $B_{1}^{H}(t),\dots,B_{p}^{H}(t)$
of a standard $d$-dimensional fractional Brownian motion $B^{H}(t)$.
Throughout this paper
\[
p^{*}:=p/(p-1)\,
\]
 will stand for the conjugate to $p>1$. Our next and main goal is to investigate large
deviations for intersection local time $\alpha^{H}(\cdot)$
of $B_{1}^{H}(t),\cdots,B_{p}^{H}(t)$,
which is a random measure on $(\R^{+})^{p}$ given heuristically by
\[
\alpha^{H}(A)=\int_{A}\,\prod_{j=1}^{p-1}\delta_{0}
\big(B_{j}^{H}(s_{j})-B_{j+1}^{H}(s_{j+1})\big)\, ds_{1}\cdots ds_{p},
\quad A\subset(\R^{+})^{p}.\]

Quantities measuring the amount of self-intersection of a random walk,
or of mutual intersection of
several independent random walks, have been studied
intensively for more than twenty years, see e.g.
\cite{DV}, \cite{LG}, \cite{La}, \cite{MR}, \cite{HK}, \cite{C},
\cite{Ch}.
This research is motivated by the role these quantities
play in quantum field theory, see e.g. \cite{FFS}, in
our understanding of self-avoiding walks and polymer
models, see e.g.\cite{ MS}, \cite{Hollander},
or in the analysis of stochastic
processes in random environments, see e.g.
\cite{HKM} \cite{GKS}, \cite{AC}, \cite{FMW}. In the
latter models dependence between a moving particle and a
random environment frequently comes from the particle's ability
to revisit sites with an attractive  (in some sense) environment. Consequently,  measures of self-intersection quantify the degree of dependence
between movement and environment. Typically, in high
dimensions, this dependence gets weaker, as
the movements become more transient and self-intersections less likely.
Investigation of large deviations for intersection local times
is closely related to asymptotics of the partition functions
in above models.

There are two equivalent ways to construct $\alpha^{H}(A)$ rigorously.
In the first way, $\alpha^{H}(A)$ is defined as the local time at zero of
the multi-parameter process
\begin{align}\label{eq:X}
X(t_1,\cdots, t_p)=
\big(B_{1}^{H}(t_{1})-B_{2}^{H}(t_{2}),\cdots,B_{p-1}^{H}(t_{p-1})
-B_{p}^{H}(t_{p})\big)\hskip.2in (t_1,\cdots, t_p)\in(\R^+)^p
\end{align}
More precisely, consider the occupation measure
\[
\mu_{A}(B)=\int_{A}\1_{B}
\big(B_{1}^{H}(s_{1})-B_{2}^{H}(s_{2}),\cdots,B_{p-1}^{H}(s_{p-1})
-B_{p}^{H}(s_{p})\big)\, ds_{1}\cdots ds_{p},\hskip.2inB\subset\R^{d(p-1)}.\]
By Theorem \ref{7.2},  as $Hd<p^{*}$,
there is a density function $\alpha^{H}(A,\,\cdot)$ of $\mu_{A}(\cdot)$
such that if $A=[0,t_{1}]\times\cdots\times[0,t_{p}]$,
then $\alpha^{H}([0,t_{1}]\times\cdots\times[0,t_{p}],x)$
is jointly continuous in $(t_{1},\cdots,t_{p},x)$.
We define $\alpha^{H}(A):=\alpha^{H}(A,\,0)$.

For the second way of constructing $\alpha^{H}(A)$, write for any
$\epsilon>0$
\begin{align}
\alpha_{\epsilon}^{H}(A)=\int_{\R^{d}}\int_{A}\,
\prod_{j=1}^{p}p_{\epsilon}\big(B_{j}^{H}(s_{j})-x\big)\, ds_{1}
\cdots ds_{p}\, dx & ,\label{eq:a_ep}
\end{align}
where $p_{\epsilon}$ are probability densities approximating $\delta_{0}$
as $\epsilon\to0$. Notice that
$$
\begin{aligned}
\alpha_{\epsilon}^{H}(A)&=\int_Ah_\epsilon
\big(B_{1}^{H}(s_{1})-B_{2}^{H}(s_{2}),\cdots,B_{p-1}^{H}(s_{p-1})
-B_{p}^{H}(s_{p})\big)\, ds_{1}\cdots ds_{p}\\
&=\int_{\R^{d(p-1)}}h_\epsilon (x)\alpha^{H}(A, x)dx
\end{aligned}
$$
where
$$
h_\epsilon (x_1,\cdots, x_{p-1})=\int_{\R^d}p_\epsilon (-x)\prod_{j=1}^{p-1}
p_\epsilon\Big(\sum_{k=j}^{p-1}x_k-x\Big)
$$ is an probability density on $\R^{d(p-1)}$ approaching
$\delta_0(x_1,\cdots, x_{p-1})$ as $\epsilon\to 0^+$.

By the continuity of $\alpha^{H}(A, x)$,
$\lim_{\epsilon\to0^{+}}\alpha_{\epsilon}^{H}(A)=\alpha^H(A)$
almost surely. Applying Proposition \ref{3.1.1} to the Gaussian field given
in (\ref{eq:X}), the convergence is also in ${\cal L}^m$ for all positive
$m$. This way
of constructing $\alpha^{H}(A)$ justifies the
symbolic notation
$$
\alpha^{H}(A)=\int_{\R^{d}}\int_{A}\,\prod_{j=1}^{p}\delta_0
\big(B_{j}^{H}(s_{j})-x\big)\, ds_{1}\cdots ds_{p}\, dx.
$$

In the special
case $p=2$
and $Hd<2$, Nualart and Ortiz-Latorre \cite{NO}
proved that $\alpha_{\epsilon}^{H}\big([0,t_{1}]\times[0,t_{2}]\big)$
converges in ${\cal L}^{2}$ as $\epsilon\to0^{+}$, with
\begin{equation}
p_{\epsilon}(x)=(2\epsilon\pi)^{-d/2} \exp\{-\vert x\vert^{2}/2\epsilon\}.
\label{eq:p_ep}
\end{equation}

For the Riemann-Liouville process $W^{H}(t)$ an analogous construction
of the intersection local time
$$
\begin{aligned}
\tilde{\alpha}^{H}(A)
&=\int_{A}\,\prod_{j=1}^{p-1}\delta_{0}
\big(W_{j}^{H}(s_{j})-W_{j+1}^{H}(s_{j+1})\big)\, ds_{1}\cdots ds_{p} \\
&=\int_{\R^{d}}\int_{A}\,\prod_{j=1}^{p}\delta_0
\big(W_{j}^{H}(s_{j})-x\big)\, ds_{1}\cdots ds_{p}\, dx,
\quad A\subset(\R^{+})^{p}
\end{aligned}
$$
can be
done under the same condition $Hd<p^*$.

By the self-similarity of $B^H(t)$ and $W^H(t)$,
for any $t>0$
\begin{align}
\alpha^{H}\big([0,t]^{p}\big)\stackrel{d}{=}t^{p-Hd(p-1)}
\alpha^{H}\big([0,1]^{p}\big) &
\label{aB-ss}
\end{align}
and
\begin{equation}
\tilde{\alpha}^{H}\big([0,t]^{p}\big)\stackrel{d}{=}t^{p-Hd(p-1)}
\tilde{\alpha}^{H}\big([0,1]^{p}\big).\label{eq:aW-ss}
\end{equation}

Finally, we would like to discuss this research in a more general context of Markovian versus non-Markovian structures. Naturally, most of the existing results on large deviation for (intersection) local time have been obtained
for Markov processes such as Brownian motions, L\'evy stable processes, general L\'evy processes, and random walks. The underlying Markovian structure has been essential for the methods in these studies; see Chen \cite{Ch} for references and a systematical account of such works. Departures from Markovian models are often driven by the underlying physics to match  the required level of dependence (memory) and smoothness/roughness of sample paths. Fractional Brownian motion and Riemann–Liouville processes are the most natural candidates as extensions of Brownian motion into the non-Markovian world. They offer the existence of the intersection local time for any number $p$ of processes in any dimension $d$ as long as $H$ is sufficiently small. Therefore, they may help scientists to build more realistic and robust models while posing serious challenge to mathematicians due to the
non-Markovian nature.

In this paper, we mainly use Gaussian techniques motivated from the 
study of continuity properties of local time, 
and more generally, from theory of Gaussian processes. It
 is also helpful to see connections between 
small ball probability estimates and tail behavior of the local time. 
Indeed, large value of the local time at zero means that the process 
stayed for a long time in a small neighborhood of zero. By this analogy, 
Propositions \ref{3.1.1} and \ref{3.1.2} can be motivated by the 
corresponding results for small balls (see comments preceding 
these propositions in Section \ref{sub:Comparison}).

\section{Main results} \label{s:main}

\begin{theorem} \label{th:1}
Let $B^{H}(t)$ be a standard $d$-dimensional
fractional Brownian motion with index $H$ such that $Hd<1$. Then
the limit
\begin{align}
\lim_{a\to\infty}a^{-1/(Hd)}\log\P\{L_{1}^{0}(B^{H})\ge a\}=-\theta(H,d)
\label{eq:limLB}
\end{align}
 exists and $\theta(H,d)$ satisfies the following bounds
 \begin{align}
\left(\pi c_{H}^{2}/H\right)^{1/(2H)}\theta_{0}(Hd)\le\theta(H,d)\le(2\pi)^{1/(2H)}\theta_{0}(Hd) & ,\label{eq:theta}
\end{align}
 where $c_{H}$ is given by (\ref{eq:c}) and
 \begin{equation}
\theta_{0}(\kappa)=\kappa \left(\frac{(1-\kappa)^{1-\kappa}}{\Gamma(1-\kappa)}\right)^{1/\kappa}\,.\label{eq:theta0}
\end{equation}
 \end{theorem}

Notice that in the classical case of one-dimensional Brownian motion,
\eqref{eq:theta} becomes the equality. The fact that the lower bound is less than or equal to the upper bound in \eqref{eq:theta} is equivalent to $c_H^2 \le 2H$, which can also be seen directly. Indeed, from (3.16)
\begin{align}\label{bounds}
\frac{c_H^2}{2H} =\mathrm{Var}(B^H(1) | B^H(s), s\le 0) \le \mathrm{Var}(B^H(1)) =1.
\end{align}
The equality only holds for a Brownian motion, i.e., $H=1/2$.

\begin{theorem} \label{th:2}
Let $W^{H}(t)$ be a $d$-dimensional
Riemann–Liouville process as in \eqref{eq:W} such that $Hd<1$. Then
the limit
\begin{align}
\lim_{a\to\infty}a^{-1/(Hd)}\log\P\{L_{1}^{0}(W^{H})\ge a\}=- & \tilde{\theta}(H,d),\label{eq:limLW}
\end{align}
 exists with
 \begin{equation}
\tilde{\theta}(H,d)=\left(c_{H}\right)^{-1/H}\theta(H,d)\,,\label{eq:tilde-theta}
\end{equation}
 where $\theta(H,d)$ is as in Theorem \ref{th:1} and $c_{H}$ is
given by (\ref{eq:c}).

\end{theorem}

\begin{theorem}\label{th:4}
Let $\tilde{\alpha}^{H}(\cdot)$ be
the intersection local time of $p$-independent $d$-dimensional Riemann–Liouville
process $W_{1}^{H}(t),\cdots,W_{p}^{H}(t)$, where $Hd<p^{*}.$ Then
the limit
\begin{align}
\lim_{a\to\infty}a^{-p^{*}/(Hdp)}\log\P\big\{\tilde{\alpha}^{H}
\big([0,1]^{p}\big)\ge a\big\}=-\widetilde{K}(H,d,p)
\label{eq:lim-aW-tilde}
\end{align}
 exists and $\tilde{K}(H,d,p)$ satisfies the following bounds
\begin{align} \label{eq:K-tilde-K}
& p\frac{Hd}{p^*}  \Big(1-\frac{Hd}{p^*}\Big)^{1-\frac{p^*}{Hd}}
\Big(\frac{\pi}{H}\Big)^{\frac{1}{2H} } p^{\frac{p^*}{2Hp}}
\Gamma\Big(1-\frac{Hd}{p^*}\Big)^{-{\frac{p^*}{Hd}}}    \le\widetilde{K}(H,d,p)\\
&  \le p \frac{Hd}{p^*} \Big(1-\frac{Hd}{p^*}\Big)^{1-\frac{p^*}{Hd} }
\left(\frac{2\pi}{c_H^2p^*}\right)^{\frac{1}{2H}}
\bigg(\int_0^\infty
\big(1+t^{2H}\big)^{-d/2}e^{-t}dt\bigg)^{-\frac{p^*}{Hd}}
\nonumber
\end{align}
where $c_{H}$ is given by (\ref{eq:c}).

\end{theorem}

There is a direct way to show that the lower bound is less than or equal to
the upper bound in (\ref{eq:K-tilde-K}).
Observe that by H\"older inequality,
$
1+t^{2H}\ge p^{1/p}(p^*)^{1/p^*} t^{2H/p^*}
$
which leads to
$$
\int_0^\infty \big(1+t^{2H}\big)^{-d/2}e^{-t}dt\le p^{-d/(2p)}(p^*)^{-d/(2p^*)}
\Gamma(1-Hd/p^*).
$$
After cancellation on both sides of (\ref{eq:K-tilde-K}), the problem is
then reduced to examining the relation $c_H^2\le 2H$, which is given in
(\ref{bounds}).

\begin{theorem}\label{th:3}
Let $\alpha^{H}(\cdot)$ be the intersection
local time of $p$-independent standard $d$-dimensional fractional
Brownian motions $B_{1}^{H}(t),\cdots,B_{p}^{H}(t)$, where $Hd<p^{*}.$
Then the limit
\begin{align}
\lim_{a\to\infty}a^{-p^{*}/(Hdp)}\log
\P\big\{\alpha^{H}\big([0,1]^{p}\big)\ge a\big\}
=-K(H,d,p)\label{eq:lim-aB}
\end{align}
 exists with
\begin{align}\label{eq:K-bound}
K(H,d,p)=c_H^{1/H}\tilde{K}(H,d,p).
\end{align}
\end{theorem}

Our results seem to be closely related to the large deviations of the
self-intersection local times heuristically
written as
$$
\beta^H\big([0,t]_<^p\big)
=\int_{[0, t]_<^p}\prod_{j=1}^{p-1}\delta_0\big(B^H(s_j)-B^H(s_{j+1})\big)
ds_1\cdots ds_p
$$
where
$$
[0,t]_<^p=\big\{(s_1,\cdots, s_p)\in [0, t]^p;\hskip.1in
s_1<\cdots <s_p\big\}.
$$
In the case when $Hd<1$, we can rewrite
$$
\beta^H\big([0,t]_<^p\big)={1\over p!}\int_{\R^d}\big[L_t^x(B^H)\big]^p \, dx.
$$
To see the connection between $\alpha^H$ and $\beta^H$,
notice that by Holder inequality and arithmetic and geometric mean inequality,
$$
\left( \alpha^H\big([0, 1]^p\big) \right)^{1/p}=\left( \int_{\R^d}
\prod_{j=1}^p L_1^x(B_j^H)dx \right)^{1/p}
\le{1\over p}\sum_{j=1}^p\bigg(\int_{\R^d}\big[L_1^x(B_j^H)\big]^p \, dx\bigg)^{1/p}.
$$
Thus, for any $\theta>0$
$$
\E\exp\Big\{\theta a^{p^*-Hd\over Hdp}\Big(
\alpha^H \big([0, 1]^p\big)\Big)^{1/p}\Big\}
\le \Bigg[\E\exp\bigg\{\theta p^{-1} a^{p^*-Hd\over Hdp}
\bigg(\int_{\R^d}\big[L_1^x(B^H)\big]^p \, dx\bigg)^{1/p}\bigg\}\Bigg]^p.
$$

On the other hand, by Theorem \ref{th:3} and Varadhan's integral lemma,
$$
\begin{aligned}
&\lim_{a\to\infty}a^{-p^*/(Hdp)}\log \E\exp\Big\{\theta p^{-1}
a^{p^*-Hd\over Hdp}\Big(
\alpha^H \big([0, 1]^p\big)\Big)^{1/p}\Big\}\\
&=\sup_{\lambda>0}\Big\{{\theta p^{-1}}\lambda^{1/p}-K(H,d,p)
\lambda^{p^*/Hdp}\Big\}\\
&=({Hd/ (p^*K(H, d, p))})^{Hd/(p^*-Hd)}(1-Hd/ p^*)
({\theta /p})^{p^*/( p^*-Hd)}.
\end{aligned}
$$
Consequently,
\begin{align}
&\liminf_{a\to\infty}a^{-p^*/(Hdp)}\log\E\exp\bigg\{\theta
a^{p^*-Hd\over Hdp}\bigg(\int_{\R^d}\big[L_1^x(B^H)\big]^p \, dx\bigg)^{1/p}\bigg\} \label{e:2-11} \\
&\ge p^{-1}
({Hd/ (p^*K(H, d, p))})^{Hd/(p^*-Hd)}(1-Hd/ p^*)
({\theta/p })^{p^*/( p^*-Hd)}. \nonumber
\end{align}

If this can be strengthened into equality with limits, then
%$$ \begin{aligned}
%&\lim_{a\to\infty}a^{-p^*/Hdp}\log\E\exp\bigg\{\theta
%a^{p^*-Hd\over Hdp}\bigg(\int_{\R^d}\big[L_1^x(B^H)\big]^pdx\bigg)^{1/p}\bigg\}\\
%&= {1\over p}\Big({Hd\over p^*K(H, d, p)}\Big)^{Hd\over p^*-Hd}
%\Big(1-{Hd\over p^*}\Big) (p\theta)^{p^*\over p^*-Hd} \end{aligned} $$
by G\"artner-Ellis theorem, for any $\lambda>0$,
$$
\begin{aligned}
&\lim_{a\to\infty}a^{-p^*/(Hdp)}\log\P\bigg\{\bigg(
\int_{\R^d}\big[L_1^x(B^H)\big]^p \, dx\bigg)^{1/p}
\ge \lambda a^{1/p}\bigg\}\\
&=-\sup_{\theta>0}\bigg\{\lambda\theta -
p^{-1} ({Hd/( p^*K(H, d, p))})^{Hd\over p^*-Hd}
(1-{Hd/ p^*})
(\theta/p)^{p^*\over p^*-Hd}\bigg\}\\
&=-p^{-1} K(H,d, p)\lambda^{p^*/(Hd)}
\end{aligned}
$$
In particular,
\begin{align}\label{conj-1}
\lim_{a\to\infty}a^{-p^*/(Hdp)}\log\P\bigg\{
\int_{\R^d}\big[L_1^x(B^H)\big]^p \, dx
\ge a\bigg\}
=-p^{-1}K(H,d, p).
\end{align}

The conjecture (\ref{conj-1}) is partially supported by a recent result of Hu, Nualart and Song (Theorem 1, \cite{HNS}) which states that
when $Hd<1$ and $p=2$
$$
\E\bigg\{\int_{\R^d}\big[L_1^x(B^H)\big]^2dx\bigg\}^n\le C^n (n!)^{Hd}
\hskip.2in n=1,2,\cdots
$$
for some $C>0$. Indeed, a standard application of Chebyshev inequality and Stirling formula leads to the upper bound of the form
$$
\lim\sup_{a\to\infty}a^{-1/(Hd)}\log\P\bigg\{
\int_{\R^d}\big[L_1^x(B^H)\big]^2dx
\ge \lambda a\bigg\}\le -l,
$$
where $l$ is a positive constant.
This rate of decay of tail probabilities is sharp by comparing it with  \eqref{e:2-11} for $p=2$.

In the case $Hd\ge 1$, $\beta^H\big([0,t]_<^p\big)$ can not be properly defined.
On the other hand,  this problem can be fixed in some cases
by renormalization. For simplicity we consider the case $p=2$. Hu and
Nualart prove (Theorem 1, \cite{HN})
that for $1\le Hd<3/2$,  the renormalized self-intersection
local time formally given as
\[
\begin{aligned}
\gamma^{H}\big([0, t]_<^2\big) & =\int\!\!\int_{\{0\le r<s\le t\}}
\delta_{0}\big(B^{H}(r)-B^{H}(s)\big) \, drds\\
 &  \quad  -\E\int\!\!\int_{\{0\le r<s\le t\}}\delta_{0}
\big(B^{H}(r)-B^{H}(s)\big) \, drds
\end{aligned}
\]
exists with the scaling property
\begin{align}
\gamma^{H}\big([0,t]_<^2\big)\stackrel{d}{=}t^{2-Hd}
\gamma^{H}\big([0,1]_<^2\big)\label{1-16}
\end{align}
We also point that an earlier work by Rosen (\cite{Rosen})
in the special case $d=2$.

Based on a similar but more heuristic reasoning, it seems plausible
to expect that
\begin{align}\label{conj-2}
\lim_{a\to\infty}a^{-1/(Hd)}\log\P\Big\{\gamma^H\big([0,1]_<^2\big)\ge a\Big\}
=-2^{(Hd)^{-1}-1}K(H,d, 2)
\end{align}
We refer the interested reader  to Theorem 4,  \cite{HNS} for
some exponential integrabilities established  by Hu, Nualart
and Song based on Clark-Ocone's formula.
We leave these problems to the future investigation.

Our large deviations estimates can be applied to
obtain the law of the iterated logarithm.

\begin{theorem}\label{1-17}
When $Hd<1$,
\begin{align}
%\limsup_{t\to\infty}\frac{L_{t}^{0}(B^{H})}{t^{1-Hd}(\log\log t)^{Hd}}
%=\theta(H,d)^{-Hd}\hskip.2ina.s.\label{1-18}
\limsup_{t\to\infty} t^{-(1-Hd)}(\log\log t)^{-Hd} L_{t}^{0}(B^{H})
=\theta(H,d)^{-Hd}\hskip.2ina.s.\label{1-18}
\end{align}
\begin{align}
\limsup_{t\to\infty} t^{-(1-Hd)}(\log\log t)^{-Hd} L_{t}^{0}(W^{H})
=\tilde{\theta}(H,d)^{-Hd}\hskip.2ina.s.\label{1-19}
\end{align}

When $Hd<p^{*}$,
\begin{align}
%\limsup_{t\to\infty}\frac{\alpha^{H}\big([0,t]^{p}\big)}{t^{p-Hd(p-1)}(\log\log t)^{Hd(p-1)}}
%=K(H, d, p)^{-Hd(p-1)}\hskip.2ina.s.
\limsup_{t\to\infty} t^{-p(1-Hd/p^*)}(\log\log t)^{-Hd(p-1)} \alpha^{H}\big([0,t]^{p}\big)
=K(H, d, p)^{-Hd(p-1)}\hskip.2ina.s.
\label{1-20}
\end{align}
\begin{align}
%\limsup_{t\to\infty}\frac{\widetilde{\alpha}^{H}\big([0,t]^{p}\big)}
%{t^{p-Hd(p-1)}(\log\log t)^{Hd(p-1)}}
\limsup_{t\to\infty} t^{-p(1-Hd/p^*)}(\log\log t)^{-Hd(p-1)} \widetilde{\alpha}^{H}\big([0,t]^{p}\big)
=\widetilde{K}(H, d, p)^{-Hd(p-1)}
\hskip.2ina.s.\label{1-21}
\end{align}
 \end{theorem}

Even with the large deviations stated in Theorem 
\ref{th:1}---\ref{th:3},
the proof of Theorem  \ref{1-17} appears to be highly non-trivial
due to long-range dependency of the model. Here we mention
some previous results given in Baraka and Mountford (\cite{BM});
Baraka, Mountford and Xiao (\cite{BMX}).
Using the large deviation estimate similar to 
(\ref{eq:limLB}), Baraka, Mountford and Xiao were able
to establish some laws of the iterated logarithm 
 which describe
the short term behaviors (as $t\to 0^+$) of the local times 
of fractional Brownian motions. As pointed out by Baraka and Mountford
(p.163, \cite{BM}), their method does not lead to the
laws of the iterated logarithm of large time given in Theorem  \ref{1-17}.

Theorem \ref{1-17} will be proved in section \ref{s:LIL}.
The proof of the lower bound appears to be highly non-trivial
due to long-range dependency of the model.
The approach
relies on a quantified use of Cameron-Martin formula.

Since all main theorems stated in this section have been known in the
classic case $H=1/2$ (see, e.g., \cite{C} and \cite{CL}), we
assume $H\not =1/2$ in the remaining of the paper.

\section{Basic Tools} \label{s:aux}

In this section we provide some basic results that will be used in our proofs. We state them separately for a convenient reference.

\subsection{Comparison of local times \label{sub:Comparison}}

We will give general comparison results for local times for Gaussian processes.
They are based on the standard Fourier analytic approach but go far beyond, motivated mainly by similar small deviation estimates.
We start with an outline of the analytic method typically used in the study of local times for Gaussian processes,
in particular on its the moments, see Berman \cite{Berman} and
Xiao \cite{Xiao-1}.

For a fixed sample function and
fixed time $t>0$, the Fourier transform on space variable $x \in \R^d$ is the function of $\lambda \in \R^d$,
$$
\int_{\R^d} e^{i \lambda \cdot x } L(t,x) dx = \int_0^t  e^{i \lambda \cdot X(s) } ds.
$$
Thus the local time $L(t,x)$ can be expressed as the inverse Fourier transform:
$$
L(t, x) = {1 \over (2\pi)^d} \int_{\R^d} e^{-i \lambda \cdot x } \int_0^t  e^{i \lambda \cdot X(s) } \, ds d\lambda.
$$
The m-th power of $L(t,x)$ is
$$
L(t, x)^m = {1 \over (2\pi)^{md} } \int_{\R^{md}} e^{-i x \cdot \sum_{k=1}^m \lambda_k }
\int_{[0,t]^m} \exp\big(i \sum_{k=1}^m \lambda_k \cdot X(s_k)\big) \, ds_1\cdots ds_m d\lambda_1\cdots d\lambda_m.
$$
Take the expected value under the sign of integration: the second exponential
in the above integral is replaced by the joint characteristic function of
$X(s_1), \cdots , X(s_m)$.
In the Gaussian case, we obtain
\begin{align*}
 & \E L(t, x)^m \\
 &= {1 \over (2\pi)^{md} }
 \int_{\R^{md}} e^{-i x \cdot \sum_{k=1}^m \lambda_k }
\int_{[0,t]^m}  \exp\Big(-\frac{1}{2}\mathrm{Var}\big(\sum_{k=1}^m \lambda_k \cdot X(s_k)\big)\Big) \, ds_1\cdots ds_m d\lambda_1\cdots d\lambda_m.
\end{align*}
Interchanging integration and applying the characteristic function inversion formula, we can get more
explicit (but somewhat less useful) expression in terms of
integration associated with
$ \det ( \E X(s_i) X(s_j) )^{-1/2}$.
Estimates of the moments of local time $L(t,x)$ thus
depend on the rate of decrease to $0$ of
$\det(\E X(s_i) X(s_j))$  as $s_j - s_{j-1} \to 0$ for some $j$.
Here in our approach, we have to make proper
adjustment by approximating $L(t,x)$.

Consider now a random fields $X({\bf t})$ taking values
in $\R^{d}$, where $\mathbf{t}=(t_{1},\ldots,t_{p})\in (\R^+)^p$.
For a fixed Borel set $A\subset (\R^+)^p$,
recall that the local time formally given
as
\begin{equation}
L_X(A,x)=\int_{A}\delta_{x}\big(X(\mathbf{s})\big)\,
d{\mathbf{s}} \label{eq:L-A}
\end{equation}
is defined as the density  of the occupation measure
$$
\mu_A(B)=\int_A1_B\big(X({\bf s})\big)d{\bf s}\hskip.2in B\subset \R^d
$$
if $\mu_A(\cdot)$ is absolutely continuous with respect to the Lebesgue
measure on $\R^d$.

Given a non-degenerate Gaussian probability density $h(x)$ on $\R^d$ and
$\epsilon>0$,
the function $h_\epsilon(x)=\epsilon^{-d/2}h(\epsilon^{-1/2}x)$ is also a
probability density. Define the smoothed local time
\begin{equation}\label{eq:L-ep}
L_X(A,x,\epsilon)=\int_{A}h_{\epsilon}\big(X({\bf s})-x\big)\,
d{\bf s}.
\end{equation}

Our first proposition provides moment comparison \eqref{eq:comp1} which can be viewed as analogy of
Anderson's inequality in the small ball analog: For independent Gaussian vectors $X$, $Y$, $X$ symmetric,
$$
\P(\|X+Y\| \le \epsilon) \le \P(\|X\| \le \epsilon).
$$
See Li and Shao \cite{LS} for various application of this useful inequality.

\begin{proposition}\label{3.1.1}
Let $A\subset (\R^+)^p$ be a fixed bounded Borel set.
Let $X({\bf t})$ ($\mathbf{t}=(t_{1},\ldots,t_{p})
\in(\R^{+})^{p}$) be a zero-mean $\R^d$-valued Gaussian random field with
the local time $L_{X}(A,x)$  continuous in $x\in\R^{d}$. Assume that
for every  $m=1,2,\dots$
\begin{align}\label{3.1.2}
\int_{A^m}d{\bf s}_1\cdots d{\bf s}_m
 \int_{(\R^d)^m}d\lambda_1\cdots d\lambda_m\exp\bigg\{-{1\over 2}
\Var\Big(\sum_{k=1}^m
\lambda_k\cdot X({\bf s}_k)\Big)\bigg\}<\infty\,.
\end{align}
 Then $L_X(A, 0)\in {\cal L}^m$ (i.e., finite m-th moment), with
\begin{align}\label{3.1.3}
\E L_X(A, 0)^m&={1\over (2\pi)^{md}}\int_{A^m}d{\bf s}_1\cdots d{\bf s}_m
 \int_{(\R^d)^m}d\lambda_1\cdots d\lambda_m \\
&\times\exp\bigg\{-{1\over 2}
\Var\Big(\sum_{k=1}^m
\lambda_k\cdot X({\bf s}_k)\Big)\bigg\}\nonumber
\end{align}
and
\begin{align}\label{3.1.4}
\lim_{\epsilon\to 0^+}\E\big\vert L_X(A,0,\epsilon)-L_X(A,0)\big\vert^m=0.
\end{align}

If $Y({\bf t})$ ($\mathbf{t}=(t_{1},\ldots,t_{p})
\in(\R^{+})^{p}$) is another
zero-mean $\R^d$-valued Gaussian random field independent of $X({\bf t})$
such that the local time $L_{X+Y}(A, x)$
of $X({\bf t})+Y({\bf t})$ is continuous in $x$, then

\begin{equation}
\E\left[L_{X+Y}(A,0)^{m}\right]\le\E\left[L_{X}(A,0)^{m}\right].
\label{eq:comp1}
\end{equation}
\end{proposition}

\proof By Fourier inversion, we have from (\ref{eq:L-ep})
%and h(x)=...
$$
L_X(A,0, \epsilon)={1\over (2\pi)^d}\int_{\R^d}d\lambda
\exp\Big\{-{\epsilon\over 2}
(\lambda\cdot \Gamma\lambda)\Big\}\int_A e^{-i\lambda\cdot X({\bf s})}d{\bf s}
$$
where $\Gamma$ is the covariance matrix of Gaussian density $h(x)$.
Using Fubini theorem,
\begin{align}\label{3.1.5}
\E L_X(A,0, \epsilon)^m
&={1\over (2\pi)^{md}}\int_{A^m}d{\bf s}_1\cdots d{\bf s}_m
 \int_{(\R^d)^m}d\lambda_1\cdots d\lambda_m\nonumber \\
&\times\exp\bigg\{-{\epsilon\over 2}
\sum_{k=1}^m\lambda_k\cdot \Gamma\lambda_k\bigg\}
\E\exp\bigg\{-i\sum_{k=1}^m
\lambda_k\cdot X({\bf s}_k)\bigg\}\\
&={1\over (2\pi)^{md}}\int_{A^m}d{\bf s}_1\cdots d{\bf s}_m
 \int_{(\R^d)^m}d\lambda_1\cdots d\lambda_m \nonumber\\
&\times\exp\bigg\{-{\epsilon\over 2}
\sum_{k=1}^m\lambda_k\cdot \Gamma\lambda_k\bigg\}
\exp\bigg\{-{1\over 2}\Var\Big(\sum_{k=1}^m
\lambda_k\cdot X({\bf s}_k)\Big)\bigg\} . \nonumber
\end{align}

By monotonic convergence theorem, the right hand side converges
to the right hand side of (\ref{3.1.3}) as $\epsilon\to 0^+$. In particular,
the family
$$
\E L_X(A,0, \epsilon)^m\hskip.2in (\epsilon>0)
$$
is bounded for $m=1,2,\cdots$. Consequently, this family is uniformly
integrable for $m=1,2,\cdots$. Therefore, (\ref{3.1.3}) and (\ref{3.1.4})
follow from the fact that  $L_X(A,0, \epsilon)$ converges to $L_X(A,0)$, which is
led by the continuity of $L_X(A, x)$.

Finally, (\ref{eq:comp1}) follows from the comparison
$$
\begin{aligned}
&\int_{A^m}d{\bf s}_1\cdots d{\bf s}_m
 \int_{(\R^d)^m}d\lambda_1\cdots d\lambda_m\exp\bigg\{-{1\over 2}
\Var\Big(\sum_{k=1}^m
\lambda\cdot\big( X({\bf s}_k)+Y({\bf s}_k)\Big)\bigg\}\\
&\le \int_{A^m}d{\bf s}_1\cdots d{\bf s}_m
 \int_{(\R^d)^m}d\lambda_1\cdots d\lambda_m\exp\bigg\{-{1\over 2}
\Var\Big(\sum_{k=1}^m
\lambda\cdot X({\bf s}_k)\Big)\bigg\}.
\end{aligned}
$$
\hfill \qed

In certain situations we can also reverse bound in \eqref{eq:comp1} as a result of the Cameron-Martin Formula. In small ball setting, this is motivated by the Chen-Li's inequality \cite{CL03} which can be used to estimate small ball probabilities under any norm
via a relatively easier $L_2$-norm estimate. See also the survey of Li and Shao \cite{LS}.
Let $X$ and $Y$ be any two centered independent Gaussian random vectors in a separable Banach space $B$ with norm $\| \cdot \|$.
We use $|\cdot|_{\mu(X)}$
to denote the inner product norm induced on $H_\mu$ by $\mu={\cal L}(X)$.
Then for any $\lambda>0$ and $\epsilon>0$,
$$
\P(\|X+ Y\| \le \epsilon) \ge \P( \| X\| \le \epsilon) \cdot
\E \exp \{ -2^{-1} |Y|^2_{\mu(X)} \},
$$
and
$$
\P(\|Y\| \le \epsilon) \ge \P( \| X\| \le \lambda \epsilon) \cdot
\E \exp \{ -2^{-1} \lambda^2 |Y|^2_{\mu(X)} \}.
$$
%In particular, for any $\lam>0$, $\e>0$ and $\delta>0$,
%$$\pr{\norm{Y} \le \e} \cdot \exp\{ -\lam^2 \delta^2/2 \}
%\ge \pr{ \norm{X} \le \lam \e} \pr{ \abs{Y}_{\mu(X)} \le \delta}.$$

Next we provide the local time counterpart of this inequality, which is crucial in our estimates. Suppose that the process $X(\mathbf{t})$,
$\mathbf{t}\in[\mathbf{0},\mathbf{T}]$, where
${\bf T}=(T_1,\dots,T_p)\in (\R_+)^p$,
can be viewed as a Gaussian random vector in a separable Banach space
$B$ such that the evaluations $x\mapsto x(\mathbf{t})$ are measurable
(say $B=C([\mathbf{0},\mathbf{T}];\R^{d})$, for concreteness). Let
$\mathcal{H}(X)$ denote the reproducing kernel Hilbert space (RKHS)
of $X(\mathbf{t})$, $\mathbf{t}\in[\mathbf{0},\mathbf{T}]$ equipped
with the norm $\|\cdot\|$. Now we will make a crucial assumption
that the independent process $Y(\mathbf{t})$,
$\mathbf{t}\in[\mathbf{0},\mathbf{T}]$
has almost all paths in $\mathcal{H}(X)$.

\begin{proposition} \label{3.1.2'}
In the above setting, under the assumptions of Proposition \ref{3.1.1},
we have
\begin{equation}
\E\left[L_{X+Y}(A,0)^{m}\right]\ge\E e^{-\frac{1}{2}\|Y\|^{2}}
\E\left[L_{X}(A,0)^{m}\right], \label{eq:comp2}
\end{equation}
for every $A\subset [{\bf 0}, {\bf T}]$ and $m\in\N$.
\end{proposition}

\proof Applying Lemma \ref{lem:A-CM}(ii),
for $g(x)=\prod_{k=1}^{m}h_{\epsilon}(x(\mathbf{s}_{k}))$,
$x\in B$, we get
\begin{align*}
\E\left[L_{X+Y}(A,0,\epsilon)^{m}\right] & =\int_{A^{m}}
d\mathbf{s}_{1}\cdots d\mathbf{s}_{m}\,\E\prod_{k=1}^{m}h_{\epsilon}
\big(X(\mathbf{s}_{k})+Y(\mathbf{s}_{k})\big)\\
& \ge  \E e^{-\frac{1}{2}\|Y\|^{2}}\int_{A^m}d\mathbf{s}_{1}\cdots
d\mathbf{s}_{m}\,\E\prod_{k=1}^{m}h_{\epsilon}\big(X(\mathbf{s}_{k}))\big)
=  \E e^{-\frac{1}{2}\|Y\|^{2}}
\E\left[L_{X}(A,0,\epsilon)^{m}\right].
\end{align*}
Applying (\ref{3.1.4}) for both processes, $X$ and $X+Y$, we get  \eqref{eq:comp2}.
\qed
\bigskip	

\subsection{RKHS of $W^{H}(t)$ and the remainder $Z^{H}(t)$ \label{sub:remainder}}

Let $H\in(0,1/2)\cup(1/2,1)$ and recall decomposition \eqref{eq:dec}:
\[
c_{H}^{-1}B^{H}(t)=W^{H}(t)+Z^{H}(t),\quad t\ge0,
\]
 where the remainder process $Z^{H}(t)$ can be written as
 \begin{equation}
Z^{H}(t)=\int_{0}^{\infty}\{(t+s)^{H-1/2}-s^{H-1/2}\}\, d\bar{B}(s),\label{eq:rem1}
\end{equation}
 with $\bar{B}(s):=B(-s)$, $s\ge0$. Clearly, $Z^{H}(t)$ is a self-similar
process with index $H$ and the processes $W^{H}(t)$ and $Z^{H}(t)$
are independent. In this section we develop a technique allowing us to treat sample paths of $Z^{H}(t)$ as, essentially, elements of the reproducing kernel Hilbert space (RKHS) of $W^{H}(t)$.

The RKHS $\mathbb{H}[0,T]$ of the the Riemann–Liouville process $\left\{ W^{H}(t)\right\} _{t\in[0,T]}$ with index $H>0$, viewed as a random element in $C[0,T]$, follows standard theory of RKHS, see \cite{LL99} and \cite{BT}.
Van der Vaart and van Zanten \cite[Lemma 10.2]{Va} proved
that 
\begin{equation}\label{eq:Va}
\mathbb{H}[0,T]=I_{0+}^{H+1/2}(L_{2}[0,T])\, ,
\end{equation}
where
\begin{equation}
I_{0+}^{\alpha}f(t)=\frac{1}{\Gamma(\alpha)}\int_{0}^{t}(t-s)^{\alpha-1}f(s)\, ds\,,\quad t\in[0,T]\label{eq:RL-op}
\end{equation}
 is the Riemann–Liouville fractional integral of order $\alpha>0$;
for $\alpha=0,$ $I_{0+}^{0}f:=f$.

\begin{proposition} \label{pro:rem}
$\left\{ Z^{H}(t)\right\} _{t\ge a}$ has $C^{\infty}$-sample paths a.s. for any $a > 0$. However, for every $T>0$
$$
\P\left( \left\{ Z^{H}(t)\right\} _{t\in[0,T]} \in \mathbb{H}[0,T] \right) =0\,.
$$
\end{proposition}

\proof Formal $n$-tuple differentiation of $Z^H(t)$ gives 
\[
\dfrac{\partial^{n}}{\partial t^{n}}Z^{H}(t)= \prod_{k=1}^n (H-\frac{2k-1}{2})\int_{0}^{\infty}(t+s)^{H-(2n+1)/2}\, d\bar{B}(s),  \quad t>0.
\]
The right hand side is a well-defined Gaussian process with locally square integrable sample paths. By consecutive integration of this process over $[a,t]$ we prove that $\left\{ Z^{H}(t)\right\} _{t\ge a}$ has $C^{(n-1)}$-sample paths, $n \ge 1$, which proves the first part of the proposition.

To prove the second part observe that 
\begin{equation}\label{eq:Z-V}
Z^H(t) = (I_{0+}^{H+1/2}V^H)(t),  \quad t \ge 0,
\end{equation}
where $V^H(t)$ is a Gaussian process given by
$$
V^H(t) = \frac{H-1/2}{\Gamma(3/2-H)}\int_{0}^{\infty} \frac{t^{-H-1/2}u^{H-1/2}}{t+u} \, d\bar{B}(u)\,,  \quad t\ge 0.
$$
Direct computation gives $\E\left[ \left(V^H(t) \right)^{2}\right]=Ct^{-1}$, where $C$ depends only on $H$. Hence $\E\|V^H\|^2_{L_2[0,T]}= \infty$ but $\E\|V^H\|_{L_1[0,T]}< \infty$. Combining the fact that $I_{0+}^{H+1/2}$ is one-to-one on $L_1[0,T]$ (see \cite[Theorem 2.4]{SKM}) with \eqref{eq:Z-V} and \eqref{eq:Va}  we get
$$
\P\left( \left\{ Z^{H}(t)\right\} _{t\in[0,T]} \in \mathbb{H}[0,T] \right)= \P\left( \left\{ V^{H}(t)\right\} _{t\in[0,T]} \in L_2[0,T] \right)  =0\,,
$$
where the last equality follows from a zero-one law and integrability of Gaussian noms.
\hfill\qed

\medskip

Direct verification whether a given function belongs to $\mathbb{H}[0,T]$ can be difficult. Therefore, we give below a simple to check sufficient condition. 
Let $AC_2^{m}[0,T]$ denote the space of functions $f$ which have continuous
derivatives up to order $m-1$ on $[0,T]$, with $f^{(m-1)}$ absolutely
continuous on $[0,T]$, and $f^{(m)} \in L_2[0,T]$, $m \in \N$.

\begin{proposition} \label{prop:RKHS0}
Let $m= \lceil H+1/2\rceil$. If  $f \in AC_{2}^{m}[0,T]$ is such that $f^{(k)}(0)=0$ for $0\le k <m$, then $f \in \mathbb{H}[0,T]$ and
\begin{equation}\label{eq:RKHSnorm}
\|f\|_{\mathbb{H}[0,T]} = k_H \|I_{0+}^{m-(H+1/2)}f^{(m)}\|_{L_2[0,T]}\,,
\end{equation}
where $k_H=\Gamma(H+1/2)^{-1}$.
\end{proposition}
\proof By our assumption $f=I_{0+}^{m}f^{(m)}$, where $f^{(m)}\in L_{2}[0,T]$. Put $g=I_{0+}^{m-(H+1/2)}f^{(m)}$. Since the operators of fractional integration $\{I_{0+}^{\alpha}:\alpha\ge0\}$ form a strongly continuous semigroup on $L_2[0,1]$ (see \cite[Theorem 2.6]{SKM}), we get that  $g \in L_2[0,T]$ and 
\begin{align*}
I_{0+}^{H+1/2}g = I_{0+}^{H+1/2} \left( I_{0+}^{m-(H+1/2)}f^{(m)} \right)=I_{0+}^{m}f^{(m)}=f \,.
\end{align*}
In view of \eqref{eq:Va}, $f \in \mathbb{H}_T$ and from \cite[Lemma 10.2]{Va} 
$$
\|f\|_{\mathbb{H}[0,T]}= k_H \|g\|_{L_{2}[0,T]} = k_H \|I_{0+}^{m-(H+1/2)}f^{(m)}\|_{L_{2}[0,T]}\,.
$$
\hfill \qed

The remainder $Z^H$ is not in $\mathbb{H}[0,T]$ by Proposition \ref{pro:rem}. The next result shows the way to circumvent this problem, which is crucial to our technique.

\begin{proposition} \label{prop:Z_a}
For any $a>0$ there is a Gaussian
process $\left\{ Z_{a}^{H}(t)\right\} _{t\ge0}$ such that
\begin{itemize}
\item [{(i)}] $Z_{a}^{H}(t)=Z^{H}(t)$ for all $t\ge a$;
\item [{(ii)}]  for any $T>0$
$$
\P\left( \left\{ Z^{H}_a(t)\right\} _{t\in[0,T]} \in \mathbb{H}[0,T] \right) =1\,.
$$
\end{itemize}
\end{proposition}

\proof First consider $H\in(0,\frac{1}{2})$, so that $m=\lceil H+1/2 \rceil=1$. Define
\[
Z_{a}^{H}(t)=\begin{cases}
At, & 0\le t\le a\\
Z^{H}(t), & t>a.\end{cases}
\]
where $A=a^{-1}Z^{H}(a)$. Since $Z_{a}^{H}(t)$ has paths in $AC_{2}^{1}[0,T]$ (see the first part of Proposition \ref{pro:rem}) and $Z_{a}^{H}(0)=0$,  (ii) holds by Proposition \ref{prop:RKHS0}.

Now we consider $H\in(\frac{1}{2},1)$, so that $m=\lceil H+1/2 \rceil=2$. Define
\[
Z_{a}^{H}(t)=\begin{cases}
B_1 t^{2}+ B_2 t^{3}, & 0\le t\le a \medskip \\
Z^{H}(t), & t>a\end{cases}
\]
where $B_1=3a^{-2}Z^{H}(a)-a^{-1}\dot{Z}^{H}(a)$, $B_2=-2a^{-3} Z^H(a)+ a^{-2}\dot{Z}^{H}(a)$, and $\dot{Z}^{H}(t):=\dfrac{\partial}{\partial t}Z^{H}(t)$. As in the previous case, part (ii) follows by Proposition \ref{prop:RKHS0}. Indeed, $Z_{a}^{H}(t)$ has paths in $AC_{2}^{2}[0,T]$, $Z^{H}_a(0)=0$, and $\dot{Z}^{H}_a(0)=0$. \hfill \qed

%:here1
The above method of modifying of $Z^H$ in a neighborhood of 0 will also be used in Section \ref{s:LIL} for other processes and the $\mathbb{H}[0,T]$-norm of such modifications will to be estimated. For this purpose the next lemma will be useful.

\begin{lemma}\label{lem:bound}
Let	$m= \lceil H+1/2\rceil$. If $f \in AC_{2}^{m}[0,T]$ and $f^{(k)}(0)=0$ for $0\le k <m$, then for every $a \in (0,T)$
$$
\|f\|_{\mathbb{H}[0,T]}^2 \le C \left\{ (T^{2m-2H}-a^{2m-2H}) \|f^{(m)}\|_{L_{\infty}[0,a]}^2 + \int_{a}^T  \left|\int_{a}^T (t-s)^{m-H-3/2}f^{(m)}(s) \, ds \right|^2  dt\right\}
$$
where $C$ depends only on $H$.
\end{lemma}
\proof
Put $\kappa= m-(H+1/2)$. In view of \eqref{eq:RKHSnorm} we get
\begin{align*}
\|f\|_{\mathbb{H}[0,T]}^2 & = k_H^2\|I_{0^+}^{\kappa}(f^{(m)}\1_{[0,a]}+f^{(m)}\1_{[a,T]})\|_{L_2[0,T]}^2 \\
&\le 2k_H^2  \|I_{0^+}^{\kappa} \1_{[0,a]}\|_{L_2[0,T]}^2 \, \|f^{(m)}\|_{L_{\infty}[0,a]}^2 +2k_H^2\|I_{0^+}^{\kappa}(f^{(m)}\1_{[a,T]})\|_{L_2[0,T]}^2 \\
& \le C (T^{2m-2H}-a^{2m-2H}) \|f^{(m)}\|_{L_{\infty}[0,a]}^2 + 2k_H^2\int_{a}^T  \left|\int_{a}^T (t-s)^{\kappa-1}f^{(m)}(s) \, ds \right|^2  dt\,.
\end{align*}
\hfill\qed

\medskip

\subsection{Technical lemmas} \label{sub:tech}

The following auxiliary results and formulas are used in the proofs
of main theorems. They are given here for a convenient reference.

\begin{lemma}\label{lem:A-CM}
Let $\mu$ be a centered Gaussian
measure in a separable Banach space $B$. Let $g:B\mapsto\R_{+}$
be a measurable function. Then
\begin{itemize}
\item [{(i)}] if $\{x\in B:g(x)\ge t\}$ is symmetric and convex for
every $t>0$, then for every $y\in B$
\[
\int_{B}g(x+y)\,\mu(dx)\le\int_{B}g(x)\,\mu(dx);
\]

\item [{(ii)}] if $g$ is symmetric ($g(-x)=g(x)$, $x\in B$), then for
every $y$ in the RKHS $\mathcal{H}_{\mu}$ of $\mu$
\[
\int_{B}g(x+y)\,\mu(dx)\ge\exp\left\{ -\frac{1}{2}\|y\|_{\mu}^{2}\right\} \int_{B}g(x)\,\mu(dx),
\]
 where $\|y\|_{\mu}$ denotes the norm in $\mathcal{H}_{\mu}$.
\end{itemize}
\end{lemma}

\proof Part \emph{(i)} follows from Anderson's inequality
\begin{align*}
\int_{B}g(x+y)\,\mu(dx)= & \int_{0}^{\infty}\mu\{x\in B:\, g(x+y)\ge t\}\, dt\\
\le & \int_{0}^{\infty}\mu\{x\in B:\, g(x)\ge t\}\, dt=\int_{B}g(x)\,\mu(dx).
\end{align*}

Part \emph{(ii)} uses Cameron-Martin formula and the convexity of
exponential function
\begin{align*}
\int_{B}g(x+y)\,\mu(dx)= & \int_{B}g(x)
\exp\left\{ \langle x,y\rangle_{\mu}
-\frac{1}{2}\|y\|_{\mu}^{2}\right\} \,\mu(dx)\\
= & \frac{1}{2}\int_{B}g(x)\exp\left\{ \langle x,y\rangle_{\mu}
-\frac{1}{2}\|y\|_{\mu}^{2}\right\} \,\mu(dx)\\
 & +\frac{1}{2}\int_{B}g(x)\exp\left\{ -\langle x,y\rangle_{\mu}
-\frac{1}{2}\|y\|_{\mu}^{2}\right\} \,\mu(dx)\\
\ge & \exp\left\{ -\frac{1}{2}\|y\|_{\mu}^{2}\right\} \int_{B}g(x)\,\mu(dx).
\end{align*}
\hfill\qed

The next lemma is well-known and goes back at least to 1950s in
equivalent forms,
see Anderson \cite[p. 42]{A58},  Berman \cite[p. 293]{B69}, and \cite[p. 71]{Berman}.
The basic fact is that conditional distribution of $X_k$ given all
the $X_i, 1 \le i <k$ is a univariate Gaussian distribution with
(conditional) mean
$\E(X_k | X_1, \dots, X_{k-1})$ and (conditional) variance
$$
\det (\mathrm{Cov}(X_{1},\dots,X_{k}))/\det(\mathrm{Cov}(X_{1},\dots,X_{k-1}))
$$
for $1\le k \le m$. 

\begin{lemma} \label{lem:Cov}
Let $(X_{1},\dots,X_{m})$ be a mean-zero
Gaussian random vector. Then
\[
\mathrm{det}(\mathrm{Cov}(X_{1},\dots,X_{m}))=\mathrm{Var}(X_{1})
\mathrm{Var}(X_{2}\,|\, X_{1})\cdots\mathrm{Var}(X_{m}\,|\, X_{m-1},\dots,X_{1}).
\]
 \end{lemma}

Let $B^{H}(t)$ be given by its moving average representation \eqref{eq:B}.
By the deconvolution formula of Pipiras and Taqqu \cite{PT} we also
have

\begin{equation}
B(t)=c_{H}^{\ast}\int_{-\infty}^{t}\left((t-s)_{+}^{1/2-H}
-(-s)_{+}^{1/2-H}\right)\, dB^{H}(s),\label{eq:PT}
\end{equation}
 where $c_{H}^{\ast}=\left\{ c_{H}\Gamma(H+1/2)\Gamma(3/2-H)\right\} ^{-1}$
and the integral with respect to $B^{H}(t)$ is well-defined in the
$L^{2}$-sense. It follows from \eqref{eq:B} and \eqref{eq:PT} that
for every $t\in\R$
\begin{equation}
{\cal F}_{t}:=\sigma\{B^{H}(s);\hskip.1in-\infty<s\le t\}=\sigma\{B(s);\hskip.1in-\infty<s\le t\},\label{eq:F_t}
\end{equation}
 where the second equality holds modulo sets of probability zero.
Then for every $s<t$
\begin{equation}
\E(B^{H}(t)\,|\,\mathcal{F}_{s})=c_{H}\int_{-\infty}^{s}\Big((t-u)^{H-\frac{1}{2}}-(-u)_{+}^{H-\frac{1}{2}}\Big)\, dB(u).\label{eq:B-proj}
\end{equation}
 If $d=1$, then for every $s<t$
\begin{align}
\Var(B^{H}(t)\,|\,\mathcal{F}_{s}) &
=\E\bigg\{\Big[B^{H}(t)-\E(B^{H}(t)\,|\,\mathcal{F}_{s})\Big]^{2}\,|
\,\mathcal{F}_{s}\bigg\}\nonumber \\
 & =\E\bigg\{\int_{s}^{t}(t-u)^{H-\frac{1}{2}}\:
dB(u)\,|\,\mathcal{F}_{s}\bigg\}\nonumber \\
 & =c_{H}^{2}\int_{s}^{t}(t-u)^{2H-1} \, du=\frac{c_{H}^{2}}
{2H}(t-s)^{2H}\,.\label{eq:var-B-proj}
\end{align}

For the reader's convenience we also quote the following lemma due
to König and Mörters, \cite[Lemma2.3]{KM02}.

\begin{lemma}\label{lem:KM}
Let $Y\ge0$ be a random variable and
let $\gamma>0$. If
\begin{equation}
\lim_{m\to\infty}\frac{1}{m}\log\left( \frac{1}{(m!)^{\gamma}}\E Y^{m}\right) =\kappa\label{eq:KM1}
\end{equation}
 for some $\kappa\in\R$, then
 \begin{equation}
\lim_{y\to\infty}\frac{1}{y^{1/\gamma}}\log\P\{Y\ge y\}=-\gamma e^{-\kappa/\gamma}.\label{eq:KM2}
\end{equation}
 \end{lemma}

\section{Large deviations for local times} \label{s:LDLT}

\subsection{Proof of Theorem \ref{th:1} -- superadditivity argument} \label{sub:Th1}

In the light of Lemma \ref{lem:KM}, it is enough to show that the
limit in (\ref{eq:KM1}) exists for $Y=L_{1}^{0}(B^{H})$ and for
$\gamma =Hd$. We will
prove it by a superadditivity argument. Let $\tau$ be an exponential
time independent of $B^{H}(t)$. We will first show that for any integer
$m,n\ge1$,
\begin{align}
\E\Big[L_{\tau}^{0}(B^{H})^{m+n}\Big]\ge{m+n \choose m}\E\Big[L_{\tau}^{0}(B^{H})^{m}\Big]\E\Big[L_{\tau}^{0}(B^{H})^{n}\Big]
\label{3-1}
\end{align}

Let $t>0$ be fixed. Notice that by Theorem \ref{7.2},
the Gaussian process $B^H(t)$ satisfies
the condition (\ref{3.1.2}) posted in Proposition \ref{3.1.1}.
By (\ref{3.1.3}), therefore,
\[
\begin{aligned}
&\E\Big[L_{t}^{0}(B^{H})^m\Big]\\
 & =
{1\over (2\pi)^{md}}\int_{[0, t]^m}ds_1\cdots ds_m
 \int_{(\R^d)^m}d\lambda_1\cdots d\lambda_m\exp\bigg\{-{1\over 2}\Var\Big(\sum_{k=1}^m
\lambda_k\cdot B^H(s_k)\Big)\bigg\}\\
 & =\frac{1}{(2\pi)^{md}}\int_{[0,t]^{m}}ds_{1}\cdots ds_{m}\bigg[\int_{\R^{m}}d\lambda_{1}\cdots d\lambda_{m}
\exp\bigg\{-{1\over 2}\Var\Big(
\sum_{k=1}^{m}\lambda_{k}B_{0}^{H}(s_{k})\Big)\bigg\}\bigg]^{d}\,
\end{aligned}
\]
where $B_{0}^{H}(t)$ is 1-dimensional fractional Brownian motion.
%$B_{0,1}^{H}(t),\cdots,B_{0,d}^{H}(t)$ are independent copies of a
%1-dimensional fractional Brownian motion $B_{0}^{H}(t)$.

By integration with respect to Gaussian measures
 \[
\begin{aligned}
&\int_{\R^{m}}d\lambda_{1}\cdots d\lambda_{m}\exp\Big\{-\frac{1}{2}\Var\Big(\sum_{k=1}^{m}\lambda_{k}B_{0}^{H}(s_{k})\Big)\Big\}\\
 & =(2\pi)^{m/2}\det\Big\{\Cov\Big(B_{0}^{H}(s_{1}),\cdots,B_{0}^{H}(s_{m})\Big)\Big\}^{-1/2}\,.
 \end{aligned}
\]
 Therefore,
 \begin{align}
\E\Big[L_{t}^{0}(B^{H})^{m}\Big]
%\nonumber\\
%& =  \frac{1}{(2\pi)^{md/2}}\int_{[0,t]^{m}}ds_{1}\cdots ds_{m}\det\Big\{\Cov\Big(B_{0}^{H}(s_{1}),\cdots,B_{0}^{H}(s_{m})\Big)\Big\}^{-d/2}\label{3-2}\\
  =  \frac{m!}{(2\pi)^{md/2}}\int_{[0,t]_{<}^{m}}ds_{1}\cdots ds_{m}\det\Big\{\Cov\Big(B_{0}^{H}(s_{1}),\cdots,B_{0}^{H}(s_{m})\Big)\Big\}^{-d/2}.\label{3-2}
 \end{align}
 In \eqref{3-2} and elsewhere, for any $A\subset\R^{+}$ and an integer
$m\ge1$, we define
\[
A_{<}^{m}\,=\big\{(s_{1},\cdots,s_{m})\in A^{m};\hskip.1ins_{1}<\cdots<s_{m}\big\}.\]
 Put
 \begin{align*}
{\cal A}(s_{1},\cdots,s_{k})=\sigma\Big\{B_{0}^{H}(s_{1}),\cdots,B_{0}^{H}(s_{k})\Big\},\hskip.2ink=1,\cdots,m,
\end{align*}
 and ${\cal A}(s_{1},\cdots,s_{k})=\{\emptyset,\Omega\}$ when $k=0$.
By Lemma \ref{lem:Cov},
\begin{align}
\E\Big[L_{t}^{0}(B^{H})^{m}\Big]=\frac{m!}{(2\pi)^{md/2}}\int_{[0,t]_{<}^{m}}ds_{1}\cdots ds_{m} \varphi_{m}(s_{1},\cdots,s_{m}),
%\,\prod_{k=1}^{m}\Var\Big(B_{0}^{H}(s_{k})|B_{0}^{H}(s_{1}),\cdots,B_{0}^{H}(s_{k-1}) \Big)^{-d/2}
\nonumber\,\label{3-3}
\end{align}
where
 \[
\varphi_{m}(s_{1},\cdots,s_{m})=\prod_{k=1}^{m}\Var\Big(B_{0}^{H}(s_{k})|B_{0}^{H}(s_{1}),\cdots,B_{0}^{H}(s_{k-1})\Big){}^{-d/2}
\]
with the convention that the first term is $\Var (B_0^H(s_1))$ for $k=1$.
We are ready to establish (\ref{3-1}). Let $m,n\ge1$ be integers. Then, for any $s_{1}<\cdots<s_{n+m}$
and $n+1\le k\le n+m$,
\[
\begin{aligned}
\Var & \Big(B_{0}^{H}(s_{k})|B_{0}^{H}(s_{1}),\cdots,B_{0}^{H}(s_{k-1})\Big)\\
 & =\Var\Big(B_{0}^{H}(s_{k})-B_{0}^H(s_{n})|B_{0}^{H}(s_{1}),\cdots,B_{0}^{H}(s_{k-1})\Big)\\
 & =\Var\Big(B_{0}^{H}(s_{k})-B_{0}^H(s_{n})|B_{0}^{H}(s_{1}),\cdots,B_{0}^{H}(s_{n}),B_{0}^{H}(s_{n+1})-B_{0}^{H}(s_{n}),\\
 & \hskip3.2in\cdots,B_{0}^{H}(s_{k-1})-B_{0}^{H}(s_{n})\Big)\\
 & \le\Var\Big(B_{0}^{H}(s_{k})-B_{0}^H(s_{n})|B_{0}^{H}(s_{n+1})-B_{0}^{H}(s_{n}),\cdots,B_{0}^{H}(s_{k-1})-B_{0}^{H}(s_{n})\Big)\\
 & =\Var\Big(B_{0}^{H}(s_{k}-s_{n})|B_{0}^{H}(s_{n+1}-s_{n}),\cdots,B_{0}^{H}(s_{k-1}-s_{n})\Big),
 \end{aligned}
\]
 where the last step follows from the stationarity of increments.
Thus
\[
\varphi_{n+m}(s_{1},\cdots,s_{n+m})\ge\varphi_{n}(s_{1},\cdots,s_{n})\varphi_{m}(s_{n+1}-s_{n},\cdots,s_{n+m}-s_{n}).
\]
 Notice that from \eqref{3-2}
 \begin{align}
\E\Big[L_{\tau}^{0}(B^{H})^{m}\Big] & =\frac{m!}{(2\pi)^{md/2}}\E\int_{[0,\tau]_{<}^{m}}ds_{1}\cdots ds_{m}\,\varphi_{m}(s_{1},\cdots,s_{m})\nonumber \\
& =\frac{m!}{(2\pi)^{md/2}}\E\int_{s_1< \cdots < s_m} 1_{s_m < \tau}ds_{1}\cdots ds_{m}\,\varphi_{m}(s_{1},\cdots,s_{m})\nonumber \\
 & =\frac{m!}{(2\pi)^{md/2}}\int_{(\R^{+})_{<}^{m}}ds_{1}\cdots ds_{m}\,\varphi_{m}(s_{1},\cdots,s_{m})e^{-s_{m}}\,.\label{3-4}
 \end{align}
 Consequently,
 \begin{align*}
\E\Big[L_{\tau}^{0}(B^{H})^{n+m}\Big] & =\frac{(n+m)!}{(2\pi)^{(n+m)d/2}}\int_{(\R^{+})_{<}^{n+m}}ds_{1}\cdots ds_{n+m}\,\varphi_{n+m}(s_{1},\cdots,s_{n+m})e^{-s_{n+m}}\\
 & \ge\frac{(n+m)!}{(2\pi)^{(n+m)d/2}}\int_{(\R^{+})_{<}^{n+m}}ds_{1}\cdots ds_{n+m}\,\varphi_{n}(s_{1},\cdots,s_{n})e^{-s_{n}}\\
 & \hskip1.3in\times\varphi_{m}(s_{n+1}-s_{n},\cdots,s_{n+m}-s_{n})e^{-(s_{n+m}-s_{n})}\\
 & =\frac{(n+m)!}{(2\pi)^{(n+m)d/2}}\int_{(\R^{+})_{<}^{n}}ds_{1}\cdots ds_{n}\,\varphi_{n}(s_{1},\cdots,s_{n})e^{-s_{n}}\\
 & \hskip1.3in\times\int_{(\R^{+})_{<}^{m}}dt_{1}\cdots dt_{m}\,\varphi_{m}(t_{1},\cdots,t_{m})e^{-t_{m}}\\
 & ={n+m \choose m}\E\Big[L_{\tau}^{0}(B^{H})^{n}\Big]\E\Big[L_{\tau}^{0}(B^{H})^{m}\Big]\,.
 \end{align*}
We proved relation \eqref{3-1} that says that the sequence
$m\mapsto\log\left( \dfrac{1}{m!}\E\Big[L_{\tau}^{0}(B^{H})^{m}\Big]\right)$
\, is super-additive. By Fekete's lemma the limit
\begin{align}
\lim_{m\to\infty}\frac{1}{m}\log\left( \frac{1}{m!}
\E\Big[L_{\tau}^{0}(B^{H})^{m}\Big]\right) =
\sup_{m\ge1}\frac{1}{m}\log \left(\frac{1}{m!}\E\Big[L_{\tau}^{0}(B^{H})^{m}\Big]\right)
=\log L\,,\label{3-5}
\end{align}
 exists, possibly as an extended number.
By the scaling property \eqref{eq:LB-ss},
\[
\begin{aligned}
\E\Big[L_{\tau}^{0}(B^{H})^{m}\Big] & =\E\Big[\tau^{(1-Hd)m}\Big]\,\E\Big[L_{1}^{0}(B^{H})^{m}\Big] 
=\Gamma\left(1+(1-Hd)m\right)\,\E\Big[L_{1}^{0}(B^{H})^{m}\Big]\,.
 \end{aligned}
\]
 From \eqref{3-5} and Stirling's formula we get
 \begin{equation}
\lim_{m\to\infty}\frac{1}{m}\log\left( \frac{1}{(m!)^{Hd}}\E\Big[L_{1}^{0}(B^{H})^{m}\Big]\right)=\log\Big\{(1-Hd)^{-(1-Hd)}L\Big\}.\label{eq:limTh1}
\end{equation}
Applying Lemma \ref{lem:KM} we establish \eqref{eq:limLB} with
\begin{equation}
\theta(H,d)=Hd(1-Hd)^{-1+1/Hd}L^{-1/Hd}\,.\label{eq:theta1}
\end{equation}
To obtain \eqref{eq:theta} and complete the proof it is enough to
show that
\begin{align}
(2\pi)^{-d/2}\Gamma(1-Hd)\le L\le\Big(H^{-1}\pi c_{H}^{2}\Big)^{-d/2}\Gamma(1-Hd)\,.\label{3-6}
\end{align}
 By (\ref{3-1})
 \[
\frac{1}{m!}\E\Big[L_{\tau}^{0}(B^{H})^{m}\Big]\ge\left\{ \E L_{\tau}^{0}(B^{H})\right\} ^{m}=\left\{(2\pi)^{-d/2}\Gamma(1-Hd)\right\} ^{m}\,,
\]
 where the equality comes from (\ref{3-4}) (for $m=1$). This proves
the lower bound in (\ref{3-6}).

To prove the upper bound, we first notice that
\begin{align}
\Var\Big(B_{0}^{H}(s_{k})\,|\, B_0^{H}(s_{1}),\cdots,B_0^{H}(s_{k-1})\Big) & \ge\Var\Big(B_{0}^{H}(s_{k})|\, B_{0}(s),s\le s_{k-1}\Big)\label{3-7}\\
 & =\frac{c_{H}^{2}}{2H}(s_{k}-s_{k-1})^{2H},\nonumber
 \end{align}
 where we used \eqref{eq:var-B-proj}. Hence the function $\varphi$
defined above satisfies, with $s_0=0$,
\[
\varphi_{m}(s_{1},\cdots,s_{m})\le\left(2H/c_{H}^{2}\right)^{md/2}\prod_{k=1}^{m}(s_{k}-s_{k-1})^{-Hd},
\]
 and by \eqref{3-4},
 \begin{align}
 \left(\pi c_{H}^{2}/H\right)^{md/2} \E\Big[L_{\tau}^{0}(B^{H})^{m}\Big] & \le m!\int_{(\R^{+})_{<}^{m}}ds_{1}\cdots ds_{m}\,\prod_{k=1}^{m}(s_{k}-s_{k-1})^{-Hd}e^{-s_{m}}\label{3-5a}\\
& =m!\left\{ \int_{0}^{\infty}t^{-Hd}e^{-t}\, dt\right\} ^{m}=m!\Gamma(1-Hd)^{m}\,.\nonumber
%\E\Big[L_{\tau}^{0}(B^{H})^{m}\Big] & \le\left(2H/c_{H}^{2}\right)^{md/2}m!\int_{(\R^{+})_{<}^{m}}ds_{1}\cdots ds_{m}\,\prod_{k=1}^{m}(s_{k}-s_{k-1})^{-Hd}e^{-s_{m}}\label{3-5a}\\
% & =\left(2H/c_{H}^{2}\right)^{md/2}m!\left\{ \int_{0}^{\infty}t^{-Hd}e^{-t}\, dt\right\} ^{m}=\left(2H/c_{H}^{2} \right)^{md/2}m!\Gamma(1-Hd)^{m}\,.\nonumber
 \end{align}
 This establishes \eqref{3-6} and completes the proof. \qed

\subsection{Proof of Theorem \ref{th:2} -- comparison argument} \label{sub:Th2}

First we note that
\begin{equation}
L_{t}^{0}\left(c_{H}^{-1}B^{H}\right)=
c_{H}^{d}L_{t}^{0}\left(B^{H}\right).\label{eq:cB}
\end{equation}
 Thus, from the decomposition \eqref{eq:dec} and \eqref{eq:comp1}
for every $m\in\N$,
\begin{equation}
c_{H}^{md}\E\left[L_{1}^{0}\left(B^{H}\right)^{m}\right]\le\E\left[L_{1}^{0}\left(W^{H}\right)^{m}\right].\label{eq:B<W}
\end{equation}
To prove a reverse inequality (up to a multiplicative constant) we
use notation \eqref{eq:L-A}. Fix $a\in(0,1)$ and let let $Z_{a}^{H}(t)$,
$t\ge0$ be the process specified in Proposition \ref{prop:Z_a} that
is also independent of $W^{H}(t),$ $t\ge0.$ We have
\[
c_{H}^{d}L_{1}^{0}\left(B^{H}\right)=L_{c_{H}^{-1}B^{H}}([0,1],0)
\ge L_{c_{H}^{-1}B^{H}}([a,1],0)=L_{W^{H}+Z_{a}^{H}}([a,1],0).\]
 Thus, by \eqref{eq:comp2} we get
\begin{align*}
c_{H}^{md}\E\left[L_{1}^{0}\left(B^{H}\right)^{m}\right]\ge & \E\left[L_{W^{H}+Z_{a}^{H}}([a,1],0)^{m}\right]\ge K_{a}\E\left[L_{W^{H}}([a,1],0)^{m}\right]\\
= & K_{a}\E\left[\left(L_{1}^{0}(W^{H})-L_{a}^{0}(W^{H})\right)^{m}\right]\\
\ge & K_{a}\left\{ \E\left[L_{1}^{0}(W^{H})^{m}\right]^{1/m}-\E\left[L_{a}^{0}(W^{H})^{m}\right]^{1/m}\right\} ^{m}\\
= & K_{a}\left(1-a^{1-Hd}\right)^{m}
\E\left[L_{1}^{0}(W^{H})^{m}\right],
\end{align*}
 where the last equality uses self-similarity \eqref{eq:LW-ss} and
$K_{a}=\E\exp\left\{ -\frac{1}{2}\|Z_{a}^{H}\|^{2}\right\} .$ Here
$\|Z_{a}^{H}\|<\infty$ a.s. is the RKHS norm associated
with $\{W^{H}(t)\}_{t\in[0,1]}$
and computed for paths of $\{Z_{a}^{H}(t)\}_{t\in[0,1]}$
. This together with \eqref{eq:B<W} yields
\[
c_{H}^{md}\E\left[L_{1}^{0}\left(B^{H}\right)^{m}\right]\le\E\left[L_{1}^{0}\left(W^{H}\right)^{m}\right]\le K_{a}^{-1}\left(1-a^{1-Hd}\right)^{-m}c_{H}^{md}\E\left[L_{1}^{0}\left(B^{H}\right)^{m}\right].\]
 Applying the limit as in \eqref{eq:limTh1} to both sides and then
passing $a\to0$ gives
\[
\lim_{m\to\infty}\frac{1}{m}\log\left( \frac{1}{(m!)^{Hd}}\E\Big[L_{1}^{0}(W^{H})^{m}\Big]\right) =\log\Big\{c_{H}^{d}(1-Hd)^{-(1-Hd)}L\Big\}.\]

Therefore, by Lemma \ref{lem:KM} the limit in \eqref{eq:limLW} exists
and $\tilde{\theta}(H,d)=c_{H}^{-1/H}\theta(H,d)$ by \eqref{eq:theta1}.
\hfill\qed

\section{Large deviations for intersection local times} \label{s:LDILT}

\subsection{Proof of Theorem \ref{th:4} -- subadditivity argument}
\label{sub:Th3}

Let $\tilde{\alpha}_{\epsilon}^{H}(A)$ be defined analogously to
\eqref{eq:a_ep} by

\[
\tilde{\alpha}_{\epsilon}^{H}(A)=\int_{\R^{d}}\int_{A}\,\prod_{j=1}^{p}p_{\epsilon}\big(W_{j}^{H}(s_{j})-x\big)\, ds_{1}\cdots ds_{p}\, dx,\]
 where $p_{\epsilon}$ is as in \eqref{eq:p_ep}. We will first prove
the subadditivity property: for every $m,n\in\N,$
\begin{align}
\E\Big[\tilde{\alpha}_{\epsilon}^{H}\big( & [0,\tau_{1}]\times\cdots\times[0,\tau_{p}]\big)^{m+n}\Big]\label{eq:suba-ae}\\
\le & {m+n \choose m}^{p}\,\E\Big[\tilde{\alpha}_{\epsilon}^{H}\big([0,\tau_{1}]\times\cdots\times[0,\tau_{p}]\big)^{m}\Big]\,
\E\Big[\tilde{\alpha}_{\epsilon}^{H}\big([0,\tau_{1}]\times\cdots\times[0,\tau_{p}]\big)^{n}\Big],\nonumber
\end{align}
where $\tau_{1},\ldots,\tau_{p}$ are iid exponential random
variables with mean 1 and independent of $W_{1}^{H}(t),\ldots,W_{p}^{H}(t).$ Indeed,
since
\[
\tilde{\alpha}_{\epsilon}^{H}\big([0,t_{1}]\times\cdots\times[0,t_{p}]\big)^{m}=\int_{(\R^{d})^{m}}dx_{1}\cdots dx_{m}\,\prod_{j=1}^{p}\prod_{k=1}^{m}\int_{0}^{t_{j}}p_{\epsilon}\big(W_{j}^{H}(s_{j,k})-x_{k}\big)\, ds_{j,k},
\]
 we can write
 \begin{equation}
\E\left[\tilde{\alpha}_{\epsilon}^{H}\big([0,\tau_{1}]\times\cdots\times[0,\tau_{p}]\big)^{m+n}\right]=\int_{(\R^{d})^{m+n}} dx_{1}\cdots dx_{m+n}\,\xi(x_{1},\dots,x_{m+n})^{p},\label{eq:a-xi}
\end{equation}
 where
 \[
\xi(x_{1},\dots,x_{m+n})=\int_{0}^{\infty}dt\, e^{-t}\int_{[0,t]^{m+n}}ds_{1}\cdots ds_{m+n}\,\E\prod_{k=1}^{m+n}p_{\epsilon}\big(W^{H}(s_{k})-x_{k}\big).
\]
Let
\[
D_{t}=\left\{ (s_{1},\ldots,s_{m+n})\in[0,t]^{m+n}:\,\max\{s_{1},\ldots,s_{m}\}\le\min\{s_{m+1},\ldots,s_{m+n}\}\right\} .
\]
There are exactly ${m+n \choose m}$ permutations $\sigma_{i}$ of
$\{1,\ldots,m+n\}$ such that $\bigcup_{i}\sigma_{i}^{-1}D_{t}=[0,t]^{m+n}$
and $\sigma_{i}^{-1}D_{t}$ are disjoint modulo sets of measure zero
(here, $\sigma(s_{1},\dots,s_{m+n}):=(s_{\sigma(1)},\dots,s_{\sigma(m+n)})$).
Therefore,
\begin{align*}
\int_{[0,t]^{m+n}} & ds_{1}\cdots ds_{m+n}\,\E\prod_{k=1}^{m+n}p_{\epsilon}\big(W^{H}(s_{k})-x_{k}\big) \\
&= \sum_{i}\int_{\sigma_{i}^{-1}D_{t}}ds_{1}\cdots ds_{m+n}\,\E\prod_{k=1}^{m+n}p_{\epsilon}\big(W^{H}(s_{k})-x_{k}\big)\\
 &= \sum_{i}\int_{D_{t}}ds_{1}\cdots ds_{m+n}\,\E\prod_{k=1}^{m+n}p_{\epsilon}\big(W^{H}(s_{k})-x_{\sigma_{i}(k)}\big),
 \end{align*}
 which gives by H\"older inequality
 \begin{align*}
\xi( & x_{1},\dots,x_{m+n})^{p}= \left\{ \sum_{i}\int_{0}^{\infty} dt\, e^{-t} \int_{D_{t}}ds_{1}\cdots ds_{m+n} \,\E\prod_{k=1}^{m+n}p_{\epsilon}\big(W^{H}(s_{k})-x_{\sigma_{i}(k)}\big)\right\}^{p}\\
 & \quad \le {m+n \choose m}^{p-1}\sum_{i} \left\{ \int_{0}^{\infty}dt\, e^{-t}\int_{D_{t}}ds_{1}\cdots ds_{m+n}\,\E\prod_{k=1}^{m+n}p_{\epsilon}\big(W^{H}(s_{k})-x_{\sigma_{i}(k)}\big)\right\}^{p}.
 \end{align*}
Substituting into \eqref{eq:a-xi} yields
\begin{align*}
& \E \Big[\tilde{\alpha}_{\epsilon}^{H}\big([0,\tau_{1}]\times\cdots\times[0,\tau_{p}]\big)^{m+n}\Big] \le {m+n \choose m}^{p-1} \sum_{i}\int_{(\R^{d})^{m+n}}dx_{1}\cdots dx_{m+n}\,\\
 &\hskip1.5in \times  \left\{ \int_{0}^{\infty}dt\, e^{-t}\int_{D_{t}}ds_{1}\cdots ds_{m+n}\,\E\prod_{k=1}^{m+n}p_{\epsilon}\big(W^{H}(s_{k})-x_{\sigma_{i}(k)}\big)\right\} ^{p}\\
 & ={m+n \choose m}^{p}\int_{(\R^{d})^{m+n}} dx_{1}\cdots dx_{m+n}\,\left\{ \int_{0}^{\infty}dt\, e^{-t}\int_{D_{t}}ds_{1}\cdots ds_{m+n}\,\E\prod_{k=1}^{m+n}p_{\epsilon}\big(W^{H}(s_{k})-x_{k}\big)\right\}^{p}.
 \end{align*}
Since the last integrand can be written as
\begin{align*}
&\left\{ \int_{0}^{\infty}dt\, e^{-t}\int_{D_{t}}ds_{1}\cdots ds_{m+n}\,\E\prod_{k=1}^{m+n}p_{\epsilon}\big(W^{H}(s_{k})-x_{k}\big)\right\}^{p} =\int_{(\R^{+})^{p}}dt_{1}\cdots dt_{p}e^{-(t_{1}+\cdots+t_{p})}\\
& \hskip1.2in  \times \int_{D_{t_{1}}\times\cdots\times D_{t_{p}}}\Big(\prod_{j=1}^{p}ds_{j,1}\cdots ds_{j,m+n}\Big)\,\E\prod_{k=1}^{m+n}\prod_{j=1}^{p}p_{\epsilon}\big(W_{j}^{H}(s_{j,k})-x_{k}\big),
\end{align*}
after integrating with respect to $x_1,\dots,x_{m+n}$ we get
\begin{align}
\E & \Big[\tilde{\alpha}_{\epsilon}^{H}\big([0,\tau_{1}]\times\cdots\times[0,\tau_{p}]\big)^{m+n}\Big]
\le {m+n \choose m}^{p} \int_{(\R^{+})^{p}}dt_{1}\cdots dt_{p}e^{-(t_{1}+\cdots+t_{p})} \label{eq:a-xi1} \\
& \times \int_{D_{t_{1}}\times\cdots\times D_{t_{p}}} \Big(\prod_{j=1}^{p}ds_{j,1}\cdots ds_{j,m+n}\Big)\, \E \prod_{k=1}^{m+n} g_{\epsilon}\big(W_1^H(s_{1,k}),\dots,W_p^H(s_{p,k})\big), \nonumber
\end{align}
where
\begin{align}
g_{\epsilon}(y_{1},\cdots,y_{p}): & =\int_{\R^{d}}\prod_{j=1}^{p}
p_{\epsilon}(y_{j}-x)\, dx\label{eq:g_e}\\
&=(2\pi \epsilon)^{-dp/2} \int_{\R^d} e^{-\left( |x|^2-2x
\cdot \overline {y}+p^{-1}\sum_{i=1}^p |y_i|^2\right)p/(2\epsilon)}\nonumber \\
 & =(2\pi\epsilon)^{-d(p-1)/2}p^{-d/2}\exp\big\{-\frac{1}{2\epsilon}
\sum_{j=1}^{p}|y_{j}-\overline{y}|^{2}\big\},\nonumber
 \end{align}
 and $\overline{y}:=p^{-1}\sum_{i=1}^{p}y_{i}$ for $y_{1},\dots,y_{p} \in \R^{d}$. Moreover,
 \begin{align}
 & \int_{D_{t_{1}} \times\cdots\times D_{t_{p}}}
\Big(\prod_{j=1}^{p}ds_{j,1}\cdots ds_{j,m+n}\Big)\,
\E \prod_{k=1}^{m+n} g_{\epsilon}\big(W_1^H(s_{1,k}),\dots,W_p^H(s_{p,k})\big) \label{eq:a-xi2}\\
&=\int_{[{\bf 0},{\bf t}]^{m}}\Big(\prod_{j=1}^{p} ds_{j,1}\cdots ds_{j,m}\Big) \int_{[{\bf 0},{\bf t}-{\bf s}^{*}]^{n}}\Big(\prod_{j=1}^{p}ds_{j,m+1}\cdots ds_{j,m+n}\Big) \nonumber \\
&\quad \times \E \prod_{k=1}^{m} g_{\epsilon}\big(W_1^H(s_{1,k}),\dots,W_p^H(s_{p,k})\big) \prod_{k=m+1}^{m+n} g_{\epsilon}\big(W_{1}^{H}(s_{1}^{*}+s_{1,k}),\cdots,W_{p}^{H}(s_{p}^{*}+s_{p,k})\big), \nonumber
\end{align}
 where
 \[
{\bf t}=(t_{1},\cdots,t_{p}),\hskip.2in{\bf s}^{*}=(s_{1}^{*},\cdots,s_{p}^{*}),
\]
 and
 \[
s_{j}^{\ast}=\max\{s_{j,k}:\,1\le k\le m\}.
\]
 Assuming that $W_{j}^{H}(t)$ are given by \eqref{eq:W} with independent
Brownian motions $B_{j}(t)$, define $\mathcal{F}_{\mathbf{s^{\ast}}}=\sigma\big\{B_{j}(u_{j}):\, u_{j}\le s_{j}^{\ast},\, j=1,\dots,p\big\}$.  Put also
\[
Y_{j}(s_{j}^{\ast},s)=\int_{s_{j}^{*}}^{s_{j}^{*}+s}(s_{j}^{*}+s-u)^{H-\frac{1}{2}}dB_{j}(u)\;\text{and}\; Z(s_{j}^{\ast},s)=\int_{0}^{s_{j}^{*}}(s_{j}^{*}+s-u)^{H-\frac{1}{2}}dB_{j}(u),
\]
so that $W_{j}(s_{j}^{\ast}+s)=Y_{j}(s_{j}^{\ast},s)+ Z_{j}(s_{j}^{\ast},s)$. The last expectation can be written as
\begin{align*}
& \E  \Big\{ \prod_{k=1}^{m} g_{\epsilon}\big(W_1^H(s_{1,k}),\dots,W_p^H(s_{p,k})\big) \\
& \times \E\Big[ \prod_{k=m+1}^{m+n} g_{\epsilon}\big(Y_{1}^{H}(s_{1}^{*}, s_{1,k}) +Z_{1}^{H}(s_{1}^{*}, s_{1,k}),\cdots,Y_{p}^{H}(s_{p}^{*},s_{p,k})+Z_{p}^{H}(s_{p}^{*},s_{p,k})\big) \big{|} \mathcal{F}_{s^{\ast}} \Big] \Big\} \\
& \le \E  \Big[ \prod_{k=1}^{m} g_{\epsilon}\big(W_1^H(s_{1,k}),\dots,W_p^H(s_{p,k})\big)\Big]\E\Big[ \prod_{k=m+1}^{m+n} g_{\epsilon}\big(Y_{1}^{H}(s_{1}^{*}, s_{1,k}),\cdots,Y_{p}^{H}(s_{p}^{*},s_{p,k})\big)  \Big]  \\
& = \E  \Big[ \prod_{k=1}^{m} g_{\epsilon}\big(W_1^H(s_{1,k}),\dots,W_p^H(s_{p,k})\big)\Big]\E\Big[ \prod_{k=m+1}^{m+n} g_{\epsilon}\big(W_{1}^{H}(s_{1,k}),\cdots,W_{p}^{H}(s_{p,k})\big)  \Big],
\end{align*}
where the inequality follows from  Lemma
\ref{lem:A-CM}(i) (see the evaluation of $g_{\epsilon}$
in \eqref{eq:g_e} and the positive quadratic form associated with it)
and the last equality follows from 
 \[
\big(Y_{1}(s_{1}^{\ast},s_{1,k}),\dots,Y_{p}(s_{p}^{\ast},s_{p,k})\big)\stackrel{d}{=}\big(W_{1}^{H}(s_{1,k}),\dots,W_{p}^{H}(s_{p,k})\big).
\]
Combining the above bound with \eqref{eq:a-xi2} and then with \eqref{eq:a-xi1} we obtain
\begin{align*}
\E & \Big[\tilde{\alpha}_{\epsilon}^{H}\big([0,\tau_{1}]\times\cdots\times[0,\tau_{p}]\big)^{m+n}\Big]
\le {m+n \choose m}^{p} \int_{(\R^{+})^{p}}dt_{1}\cdots dt_{p}e^{-(t_{1}+\cdots+t_{p})}  \\
&\times \int_{[{\bf 0},{\bf t}]^{m}}\Big(\prod_{j=1}^{p} ds_{j,1}\cdots ds_{j,m}\Big) \E \prod_{k=1}^{m} g_{\epsilon}\big(W_1^H(s_{1,k}),\dots,W_p^H(s_{p,k})\big) \\
&\times  \int_{[{\bf 0},{\bf t}-{\bf s}^{*}]^{n}}\Big(\prod_{j=1}^{p}ds_{j,m+1}\cdots ds_{j,m+n}\Big)
\E \prod_{k=m+1}^{m+n} g_{\epsilon}\big(W_{1}^{H}(s_{1,k}),\cdots,W_{p}^{H}(s_{p,k})\big)   \\
&   ={m+n \choose m}^{p}\int_{(\R_{+})^{m}} \Big(\prod_{j=1}^{p} ds_{j,1}\cdots ds_{j,m}\Big) \E \prod_{k=1}^{m} g_{\epsilon}\big(W_1^H(s_{1,k}),\dots,W_p^H(s_{p,k})\big) \\
 & \quad \times e^{-(s_{1}^{\ast}+\cdots+s_{p}^{\ast})}\, \int_{[\mathbf{s^{\ast}},\infty]^{p}}dt_{1}\cdots dt_{p}e^{-[(t_{1}-s_{1}^{\ast})+\cdots+(t_{p}-s_{p}^{\ast})]}\\
& \quad \times \int_{[\mathbf{0},\mathbf{t-s^{\ast}}]^{n}}\Big(\prod_{j=1}^{p}ds_{j,m+1}\cdots ds_{j,m+n}\Big) \E \prod_{k=m+1}^{m+n} g_{\epsilon}\big(W_{1}^{H}(s_{1,k}),\cdots,W_{p}^{H}(s_{p,k})\big) \\
 & ={m+n \choose m}^{p}\,\E\Big[\tilde{\alpha}_{\epsilon}^{H}\big([0,\tau_{1}]\times\cdots\times[0,\tau_{p}]\big)^{m}\Big]\,
\E\Big[\tilde{\alpha}_{\epsilon}^{H}\big([0,\tau_{1}]\times\cdots
\times[0,\tau_{p}]\big)^{n}\Big],
\end{align*}
 where in the last equality we use
\[
e^{-(s_{1}^{\ast}+\cdots+s_{p}^{\ast})}=\int_{(\R_{+})^{p}}
e^{-(t_{1}+\cdots+t_{p})} \prod_{k=1}^m
\1_{[\mathbf{s}^*,\mathbf{t}]}(s_{1,k},\dots,s_{p,k})\, dt_{1}\cdots dt_{p}
\]
and the definition of $g_\epsilon$ in \eqref{eq:g_e}.
The subadditivity \eqref{eq:suba-ae} is thus proved for any $\epsilon>0$.

Now we would like to take  $\epsilon\to 0^+$ on the both sides
of \eqref{eq:suba-ae}
in an attempt to establish
\begin{align}\label{eq:suba-a}
&\E\tilde{\alpha}^{H}\big([0,\tau_{1}]
\times\cdots\times[0,\tau_{p}]\big)^{m+n}\\
&\le {m+n \choose n}^{p}\,\E\tilde{\alpha}^{H}\big([0,\tau_{1}]
\times\cdots\times[0,\tau_{p}]\big)^{m}\,\E\tilde{\alpha}^{H}
\big([0,\tau_{1}]\times\cdots\times[0,\tau_{p}]\big)^{n}.\nonumber
\end{align}
To this end we need to show that for any $m\ge 1$,
$\tilde{\alpha}^{H}\big([0,\tau_{1}]\times\cdots\times[0,\tau_{p}]\big)$
is indeed in ${\cal L}^m(\Omega, {\cal A}, \P)$ and
\begin{align}\label{5.1.1}
\lim_{\epsilon\to 0^+}
\E\Big[\tilde{\alpha}_{\epsilon}^{H}\big([0,\tau_{1}]\times\cdots
\times[0,\tau_{p}]\big)^{m}\Big]
=\E\Big[\tilde{\alpha}^{H}\big([0,\tau_{1}]
\times\cdots\times[0,\tau_{p}]\big)^{m}\Big].
\end{align}
Indeed, using \eqref{eq:suba-ae} repeatedly we have that
$$
\E\Big[\tilde{\alpha}_\epsilon^{H}\big([0,\tau_{1}]
\times\cdots\times[0,\tau_{p}]\big)^{m}\Big]
\le (m!)^p\E\Big[\tilde{\alpha}_\epsilon^{H}\big([0,\tau_{1}]
\times\cdots\times[0,\tau_{p}]\big)\Big]\, .
$$
Notice that
$$
\begin{aligned}
&\E\Big[\tilde{\alpha}_\epsilon^{H}\big([0,\tau_{1}]
\times\cdots\times[0,\tau_{p}]\big)\Big]
=\int_{\R^d}\bigg[\int_0^\infty e^{-t}\E p_\epsilon \big(W^H(t)-x\big)dt\bigg]^pdx\\
&=\int_{\R^d}\bigg[\int_0^\infty e^{-t}
\int_{\R^d}p_\epsilon(y-x)p_{t^*}(y) \, dy\bigg]^p \, dx
\end{aligned}
$$
where $t^*=(2H)^{-1}t^{2H}$ and the last step
follows from the easy-to-check
fact that $W^H(t)\sim N(0, (2H)^{-1}t^{2H}I_d)$.
By Jensen inequality, the right hand side is less than or equal to
$$
\begin{aligned}
&\int_{\R^d}\int_{\R^d}p_\epsilon(y-x)\bigg[\int_0^\infty e^{-t}
\int_{\R^d}p_{t^*}(y) \, dy\bigg]^p \, dydx \\
%=\int_{\R^d}\bigg[\int_0^\infty e^{-t}\int_{\R^d}p_{t^*}(y) \, dy\bigg]^p \, dy\\
&=\int_0^\infty\!\!\cdots\!\!
\int_0^\infty dt_1\cdots dt_p
e^{-(t_1+\cdots +t_p)}\int_{\R^d}\prod_{j=1}^p p_{t_j^*}(x)dx\\
&=\Big(H/\pi\Big)^{d(p-1)/2}\int_0^\infty\!\!\cdots\!\!
\int_0^\infty e^{-(t_1+\cdots +t_p)}\Big(\sum_{j=1}^p\prod_{ 1\le k\neq j\le p} t_k^{2H}\Big)^{-d/2}dt_1\cdots dt_p
\end{aligned}
$$
where the last step follows from a routine Gaussian integration.

By arithmetic-geometric mean inequality,
$$
{1\over p}\sum_{j=1}^p\prod_{ 1\le k\neq j\le p}t_k^{2H}\ge
\prod_{j=1}^p\prod_{ 1\le k\neq j\le p}t_k^{2H/p}
=\prod_{j=1}^pt_k^{2H(p-1)/ p}.
$$
So we have
$$
\begin{aligned}
\E\Big[\tilde{\alpha}_\epsilon^{H}\big([0,\tau_{1}]
\times\cdots\times[0,\tau_{p}]\big)\Big]
&\le \Big({H/ \pi}\Big)^{d(p-1)/2}
p^{-d/2}\bigg(\int_0^\infty
t^{-{Hd(p-1)/p}}e^{-t}dt\bigg)^p\\
&=\Big({H/\pi}\Big)^{d(p-1)/2}
p^{-d/2}\Gamma(1-Hd/p^*)^p.
\end{aligned}
$$
Summarizing our computation, we obtain
\begin{align}\label{5.1.2}
(m!)^{-p}\E\Big[\tilde{\alpha}_\epsilon^{H}\big([0,\tau_{1}]
\times\cdots\times[0,\tau_{p}]\big)^m\Big]
\le \bigg(\Big({H/ \pi}\Big)^{d(p-1)/2}
p^{-d/2}\Gamma(1-Hd/p^*)^p\bigg)^m.
\end{align}
By Theorem \ref{7.2}, the process
$$
X^H(t_1,\cdots, t_p)=\Big(W_{1}^H(t_1)-W_2^H(t_2),
\cdots,W_{p-1}^H(t_{p-1})-W_p^H(t_p)\Big)
$$
satisfies the condition (\ref{3.1.2}) with $A=[{\bf 0}, {\bf t}]
=[0, t_1]\times\cdots [0, t_p]$ for any $t_1,\cdots, t_p\ge 0$
and
$$
\begin{aligned}
&\tilde{\alpha}_{\epsilon}^{H}\big([0,t_{1}]\times\cdots
\times[0,t_{p}]\big)\\
&=\int_{[{\bf 0}, {\bf t}]}
h_\epsilon\big(W_1^H(s_1)-W_2^H(s_2),\cdots,W_{p-1}^H(s_{p-1})-W_p^H(s_p)\big)
ds_1\cdots ds_p
\end{aligned}
$$
where
$$
h_\epsilon(x_1,\cdots, x_{p-1})=\int_{\R^d}p_\epsilon (-x)\prod_{j=1}^{p-1}p_1\Big(\sum_{k=j}^{p-1}
x_k-x\Big)dx
$$
is a non-degenerate normal density on $\R^{d(p-1)}$. By
Proposition \ref{3.1.1}, $\tilde{\alpha}_{\epsilon}^{H}\big([0,t_{1}]\times\cdots
\times[0,t_{p}]\big)\in{\cal L}^m(\Omega, {\cal A}, \P)$ and
\begin{align}\label{5.1.3}
\lim_{\epsilon\to 0^+}
\E\Big[\tilde{\alpha}_{\epsilon}^{H}\big([0,t_{1}]\times\cdots
\times[0,t_{p}]\big)^{m}\Big]
=\E\Big[\tilde{\alpha}^{H}\big([0,t_{1}]
\times\cdots\times[0,t_{p}]\big)^{m}\Big]\, .
\end{align}

In addition, by the representation (\ref{3.1.5}) one can see that
for any $\epsilon'<\epsilon$,
$$
\E\Big[\tilde{\alpha}_{\epsilon}^{H}\big([0,t_{1}]\times\cdots
\times[0,t_{p}]\big)^{m}\Big]
\le \E\Big[\tilde{\alpha}_{\epsilon'}^{H}\big([0,t_{1}]\times\cdots
\times[0,t_{p}]\big)^{m}\Big]\, .
$$
Thus, (\ref{5.1.1}) follows from monotonic convergence theorem and
the identities
\begin{align}\label{5.1.4}
&\E \Big[ \tilde{\alpha}_{\epsilon}^{H}\big([0,\tau_{1}]\times\cdots
\times[0,\tau_{p}]\big)^m \Big] \\
&=\int_{(\R^+)^p}e^{-(t_1+\cdots +t_p)}\E\Big[\tilde{\alpha}_{\epsilon}^{H}\big([0,t_{1}]\times\cdots
\times[0,t_{p}]\big)^{m}\Big]dt_1\cdots dt_p\nonumber
\end{align}
and
\begin{align*}%\label{5.1.5}
\E\Big[\tilde{\alpha}^{H}\big([0,\tau_{1}]\times\cdots
\times[0,\tau_{p}]\big)^m\Big]
=\int_{(\R^+)^p}e^{-(t_1+\cdots +t_p)}\E\Big[\tilde{\alpha}^{H}\big([0,t_{1}]\times\cdots
\times[0,t_{p}]\big)^{m}\Big]dt_1\cdots dt_p . \nonumber
\end{align*}
Further, by (\ref{5.1.2}) we obtain the bound
\begin{align}\label{5.1.6}
(m!)^{-p}\E\Big[\tilde{\alpha}^{H}\big([0,\tau_{1}]\times\cdots
\times[0,\tau_{p}]\big)^m\Big]
\le \bigg(\Big({H/ \pi}\Big)^{d(p-1)/2}
p^{-d/2}\Gamma(1-Hd/p^*)^p\bigg)^m \, .
\end{align}

The inequality (\ref{eq:suba-a}) implies
 that the sequence $m\mapsto\log
\left( (m!)^{-p}\E\tilde{\alpha}^{H}\big([0,\tau_{1}]
\times\cdots\times[0,\tau_{p}]\big)^{m} \right)$
\, is sub-additive.
Hence the limit
\begin{align}
\lim_{m\to\infty}\frac{1}{m}\log \left( (m!)^{-p}\E\tilde{\alpha}^{H}\big([0,\tau_{1}]\times\cdots\times[0,\tau_{p}]\big)^{m}\right) &
=c(H,d,p)\,\label{eq:lim-atau}
\end{align}
exists, possibly as an extended number. Further, by (\ref{5.1.6})
\begin{align}\label{5.1.7}
c(H,d,p)\le \log\bigg\{\Big({H/\pi}\Big)^{d(p-1)/2}
p^{-d/2}\Gamma(1-Hd/p^*)^p\bigg\}\, .
\end{align}

Now we will deduce the moments behavior
of $\tilde{\alpha}^{H}\big([0,1]^{p}\big)$.

Notice that $\tau_{*}=\min\{\tau_{1},\cdots,\tau_{p}\}$ is
an exponential time with parameter $p$.
\begin{align*}
\E\tilde{\alpha}^{H}\big( & [0,\tau_{1}]
\times\cdots\times[0,\tau_{p}]\big)^{m}
\ge\E\tilde{\alpha}^{H}\big([0,\tau_{*}]^{p}\big)^{m}\\
 & =\E\tau_{*}^{(p-Hd(p-1))m}\E\tilde{\alpha}^{H}\big([0,1]^{p}\big)^{m}\\
 & =p^{-(p-Hd(p-1))m}\Gamma\big(1+(p-Hd(p-1))m\big)
\E\tilde{\alpha}^{H}\big([0,1]^{p}\big)^{m}.
\end{align*}
By Stirling's formula,
\begin{align*}
\limsup_{m\to\infty}\frac{1}{m}\log 
\left\{ (m!)^{-Hd(p-1)}
\E\tilde{\alpha}^{H}\big([0,1]^{p}\big)^{m}\right\} 
 \le c(H,d,p)-p(1-Hd/p^*)\log(1-Hd/p^*).
\end{align*}
On the other hand, for every $t_{1},\ldots,t_{p}>0,$
\begin{align*}
\E\tilde{\alpha}_{\epsilon}^{H}
& \big([0,t_{1}]\times\cdots\times[0,t_{p}]\big)^{m}\\
= & \int_{(\R^{d})^{m}}dx_{1}\cdots dx_{m}\prod_{j=1}^{p}\int_{[0,t_{j}]^{m}}ds_{1}\cdots ds_{m}\E\prod_{k=1}^{m}p_{\epsilon}\big(W^{H}(s_{k})-x_{k}\big)\\
\le & \prod_{j=1}^{p}\left\{ \int_{(\R^{d})^{m}}dx_{1}\cdots dx_{m}\left(\int_{[0,t_{j}]^{m}}ds_{1}\cdots ds_{m}\E\prod_{k=1}^{m}p_{\epsilon}\big(W^{H}(s_{k})-x_{k}\big)\right)^{p}\right\} ^{1/p}\\
= & \prod_{j=1}^{p}\left\{ \E\tilde{\alpha}_{\epsilon}^{H}\big([0,t_{j}]^{p}\big)^{m}\right\} ^{1/p}.
\end{align*}
 Letting $\epsilon\to 0^+$,  from (\ref{5.1.3}) we get
\[
\E\tilde{\alpha}^{H}\big([0,t_{1}]\times\cdots\times[0,t_{p}]\big)^{m}
\le\prod_{j=1}^{p}\left\{ \E\tilde{\alpha}^{H}\big([0,t_{j}]^{p}\big)^{m}\right\} ^{1/p}
=\E\tilde{\alpha}^H\big([0,1]^{p}\big)^{m} \cdot \prod_{j=1}^{p}t_{j}^{m(1-Hd/p^*)} , \]
 where the last equality uses self-similarity \eqref{eq:aW-ss}.
 Hence
\begin{align}\label{5.1.7'}
\E\tilde{\alpha}^{H}\big( & [0,\tau_{1}]\times\cdots\times
[0,\tau_{p}]\big)^{m}\\
 & =\int_{(\R_{+})^{p}}dt_{1}\cdots dt_{p}\, e^{-(t_{1}+\cdots+t_{p})}
\E\tilde{\alpha}^{H}\big([0,t_{1}]\times\cdots\times[0,t_{p}]\big)^{m}\nonumber\\
 & \le\E\tilde{\alpha}^{m}\big([0,1]^{p}\big)^{m}\int_{(\R_{+})^{p}}dt_{1}\cdots dt_{p}\, e^{-(t_{1}+\cdots+t_{p})}(t_{1}\cdots t_{p})^{m(1-Hd/p^*)}\nonumber\\
 & =\E\tilde{\alpha}^{H}\big([0,1]^{p}\big)^{m}
\Gamma\Big(1+m(1-Hd/p^*) \Big)^{p}.\nonumber
 \end{align}
By Stirling's formula again,
\begin{align*}
  \liminf_{m\to\infty}\frac{1}{m}
\log\left\{ (m!)^{-Hd(p-1)}\E\tilde{\alpha}^{H}\big([0,1]^{p}\big)^{m}\right\}
 \ge c(H,d,p)-p(1-Hd/p^*)\log(1-Hd/p^*).
 \end{align*}
We have shown that
\begin{equation}
\lim_{m\to\infty}\frac{1}{m}\log\left\{ (m!)^{-Hd(p-1)}\E\tilde{\alpha}^{H}\big([0,1]^{p}\big)^{m}\right\} =C(H,d,p),\label{eq:limTh4}
\end{equation}
where by (\ref{5.1.7}),
\begin{align}\label{5.1.8}
C(H,d,p)&=c(H,d,p)-p(1-Hd/p^*)\log(1-Hd/p^*)\\
&\le\log\left\{\Big({H/ \pi}\Big)^{d(p-1)/2}
p^{-d/2}\Gamma(1-Hd/p^*)^p (1-Hd/p^*)^{-p(1-Hd/p^*)}\right\}\, .
\nonumber
\end{align}

On the other hand, let $\bar{\alpha}^H(A)$ be the intersection local time
generated by $c_H^{-1}B^H_1(t),\cdots, c_H^{-1}B^H_p(t)$. We have that
\begin{align}\label{5.1.9}
\bar{\alpha}^H(A)=c_H^{d(p-1)}\alpha^H (A), \hskip.2in A\subset (\R^+)^p.
\end{align}
In view of the decomposition (\ref{eq:dec}), by Proposition \ref{3.1.1} we have that
\begin{align}\label{5.1.10}
\E\Big[\tilde{\alpha}^H\big([0,1]^p\big)^m\Big]\ge
\E\Big[\bar{\alpha}^H\big([0,1]^p\big)^m\Big]
=c_H^{d(p-1)m}
\E\Big[\alpha^H\big([0,1]^p\big)^m\Big].
\end{align}
It follows from (\ref{5.2.6}) below that
\begin{align}\label{5.1.11}
C(H,d,p)&\ge p \log\bigg\{c_H^{d/p^*} (1-Hd/p^*)^{-(1-Hd/p^*)}
\left({p^* \over 2\pi} \right)^{\frac{d}{2p^*}}
\int_0^\infty \big(1+t^{2H}\big)^{-d/2}e^{-t}dt
\bigg\}\, .
\end{align}

Applying Lemma \ref{lem:KM} leads
the first conclusion \eqref{eq:lim-aW-tilde} of our theorem with
$$
\tilde{K}(H, d, p)=Hd(p-1)\exp\Big\{-{C(H,d,p)\over Hd(p-1)}\Big\}
$$
and therefore the bounds given in
(\ref{eq:K-tilde-K}) follows from (\ref{5.1.8})
and (\ref{5.1.11}). \hfill \qed

\subsection{Proof of Theorem \ref{th:3}  -- comparison argument} 
\label{sub:Th4}

In connection to (\ref{eq:limTh4}), we first show that
\begin{equation}
\lim_{m\to\infty}\frac{1}{m}\log\left\{ (m!)^{-Hd(p-1)}
\E\alpha^{H}\big([0,1]^{p}\big)^{m}\right\} =C(H,d,p)-d(p-1)\log c_H
,\label{eq:limTh3}
\end{equation}
The upper bound
%\begin{equation}\limsup_{m\to\infty}\frac{1}{m}
%\log\left\{ \frac{1}{(m!)^{Hd(p-1)}}
%\E\alpha^{H}\big([0,1]^{p}\big)^{m}\right\} \le C(H,d,p)
%-d(p-1)\log c_H ,\label{5.2.1} 
%\end{equation}
follows immediately from (\ref{eq:limTh4}) and the comparison (\ref{5.1.10}).
To establish the lower bound, we once again consider the intersection
local time $\bar{\alpha}^H(A)$
generated by the normalized fractional Brownian motions
$$
\bar{B}_1^H(t)=c_H^{-1}B_1^H(t),\cdots, \bar{B}_p^H(t)=c_H^{-1}B_p^H(t).
$$
For any
$\epsilon>0$, define
\begin{align}
\bar{\alpha}_{\epsilon}^{H}(A)
=\int_{\R^{d}}\int_{A}\,\prod_{j=1}^{p}p_{\epsilon}
\big(\bar{B}_{j}^{H}(s_{j})-x\big)\, ds_{1}\cdots ds_{p}\, dx & ,\label{5.2.2}
\end{align}
Let $0<\delta<1$ be a small but fixed number. Notice
$$
\begin{aligned}
\E\bar{\alpha}^H_\epsilon\big([0, 1]^p\big)^m
\ge \E\alpha^H_\epsilon\big([\delta, 1]^p\big)^m
%&=\int_{(\R^d)^m}dx_1\cdots dx_m\prod_{j=1}^p\E\int_\delta^1p_\epsilon \big(\bar{B}_j(s_k)-x_k\big)\\
=\int_{([\delta, 1]^p)^m}d{\bf s}_1\cdots d{\bf s}_m
\E\prod_{k=1}^m g_\epsilon \Big(\bar{B}_1^H(s_{1,k}),\cdots,
\bar{B}_p^H(s_{p, k})\Big)
\end{aligned}
$$
where $g_\epsilon(x_1,\cdots, x_p)$ is defined by (\ref{eq:g_e}) and
we adopt the notation ${\bf s_k}=(s_{1, k},\cdots s_{p, k})$.

Consider $\big(W_1^H(t_1),\cdots, W_p^H(t_p)\big)$
(${\bf t}=(t_1,\cdots, t_p)\in [0,1]^p$) as a Gaussian random variable
taking values in the Banach space $\otimes_{j=1}^p
C\big\{[0, 1]^p, \R^d\big\}$. Then the reproducing kernel Hilbert space
of $\Big(W_1^H(t_1),\cdots, W_p^H(t_p)\Big)$ is
$\widetilde{H}_W=\otimes_{j=1}^pH_W$.
For each $\big(f_1(t_1)\cdots, f_p(t)\big)\in\widetilde{H}_W$
$$
\big\vert\big\vert \big(f_1(t_1)\cdots, f_p(t)\big)
\big\vert\big\vert_{\widetilde{H}_W}^2
=\sum_{j=1}^p\vert\vert f_j\vert\vert_{H_W}^2
$$
where $\vert\vert\cdot \vert\vert_{H_W}$ is the reproducing kernel
Hilbert norm of $H_W$.

Let $Z_{\delta, 1}^H(t),\cdots, Z_{\delta, p}^H(t)$ be the  processes
constructed in Lemma \ref{prop:Z_a} (with $a=\delta$) by
$Z_{1}^H(t),\cdots, Z_{p}^H(t)$, respectively. For each
$({\bf s}_1,\cdots, {\bf s}_m)\in [\delta, 1]^p)^m$ by the decomposition
(\ref{eq:dec}) we have
$$
\begin{aligned}
&\E\prod_{k=1}^m g_\epsilon \Big(\bar{B}_1^H(s_{1,k}),
\cdots, \bar{B}_p^H(s_{p, k})\Big)\\
%&=\E\prod_{k=1}^m g_\epsilon \Big(W_1^H(s_{1,k})+Z_1^H(s_{1,k}),\cdots, W_p^H(s_{p, k})+Z_p^H(s_{p, k})\Big)\\
&=\E\prod_{k=1}^m g_\epsilon \Big(W_1^H(s_{1,k})+Z_{\delta, 1}^H(s_{1,k}),\cdots,
W_p^H(s_{p, k})+Z_{\delta, p}^H(s_{p, k})\Big)
\end{aligned}
$$
Fixed $({\bf s}_1,\cdots, {\bf s}_m)\in [\delta, 1]^p)^m$.
Applying Lemma \ref{lem:A-CM}(ii) to the functional $g(f_1,\cdots, f_p)$
on $\otimes_{j=1}^p
C\big\{[0, 1]^p, \R^d\big\}$ defined by
$$
g(f_1,\cdots, f_p)\equiv\prod_{k=1}^m g_\epsilon
\Big(f_1(s_{1,k}),\cdots,f_p(s_{p,k}) \Big)\hskip.2in
(f_1,\cdots, f_p)\in\otimes_{j=1}^p
C\big\{[0, 1]^p, \R^d\big\},
$$
then the right hand side is greater than
$$
\begin{aligned}
&\bigg(\E\exp\Big\{-{1\over 2}
\vert\vert Z_\delta^H\vert\vert_{H_W}^2\Big\}\bigg)^p
\E g\big(W_1^H,\cdots, W_p^H\big)\\
&=\bigg(\E\exp\Big\{-{1\over 2}
\vert\vert Z_\delta^H\vert\vert_{H_W}^2\Big\}\bigg)^p
\E\prod_{k=1}^m g_\epsilon \Big(W_1^H(s_{1,k}),
\cdots, W_p^H(s_{p, k})\Big).
\end{aligned}
$$

Summarizing our estimate, we have
$$
\begin{aligned}
&\E\bar{\alpha}^H_\epsilon\big([0, 1]^p\big)^m\\
&\ge \bigg(\E\exp\Big\{-{1\over 2}
\vert\vert Z_\delta^H\vert\vert_{H_W}^2\Big\}\bigg)^p
\int_{([\delta, 1]^p)^m}d{\bf s}_1\cdots d{\bf s}_m
\E\prod_{k=1}^m g_\epsilon \Big(W_1^H(s_{1,k}),
\cdots, W_p^H(s_{p, k})\Big)\\
&=\bigg(\E\exp\Big\{-{1\over 2}
\vert\vert Z_\delta^H\vert\vert_{H_W}^2\Big\}\bigg)^p
\E\tilde{\alpha}^H_\epsilon\big([\delta, 1]^p\big)^m.
\end{aligned}
$$
By
Proposition  \ref{3.1.1}, letting $\epsilon\to 0^+$ on both sides
yields
$$
\E\bar{\alpha}^H\big([0, 1]^p\big)^m
\ge \bigg(\E\exp\Big\{-{1\over 2}
\vert\vert Z_\delta^H\vert\vert_{H_W}^2\Big\}\bigg)^p
\E\tilde{\alpha}^H_\epsilon\big([\delta, 1]^p\big)^m.
$$

In view of (\ref{5.1.9}),
\begin{align}\label{5.2.4}
&\liminf_{m\to\infty}\frac{1}{m}\log \left( (m!)^{-Hd(p-1)}
\E\Big[\alpha^{H}\big([0,1]^{p}\big)^{m}\Big]\right) \\
&\ge -d(p-1)\log c_H+
\liminf_{m\to\infty}\frac{1}{m}\log \left( (m!)^{Hd(p-1)}
\E\Big[\tilde{\alpha}^{H}\big([\delta,1]^{p}\big)^{m}\Big] \right).\nonumber
\end{align}

To establish the lower bound for (\ref{eq:limTh3}), therefore, it remains
to show that
\begin{align}\label{5.2.5}
\liminf_{\delta\to 0^+}
\liminf_{m\to\infty}\frac{1}{m}\log \frac{1}{(m!)^{Hd(p-1)}}
\E\Big[\tilde{\alpha}^{H}\big([\delta,1]^{p}\big)^{m}\Big]
\ge C(H, d, p).
\end{align}

Write
$$
\tilde{\alpha}^H\big([0,1]^p\big)=
\tilde{\alpha}^H\big([\delta, 1]\times [0,1]^{p-1}\big)
+\tilde{\alpha}^H\big([0,\delta]\times [0,1]^{p-1}\big).
$$
By triangular inequality,
$$
\begin{aligned}
&\bigg\{\E\Big[\tilde{\alpha}^H\big([0,1]^p\big)^m\Big]\bigg\}^{1/m}\\
&\le\bigg\{\E\Big[\tilde{\alpha}^H\big([\delta, 1]
\times [0,1]^{p-1}\big)^m\Big]\bigg\}^{1/m}
+\bigg\{\E\Big[\tilde{\alpha}^H\big([0, \delta]\times
[0,1]^{p-1}\big)^m\Big]\bigg\}^{1/m}.
\end{aligned}
$$
Given $\epsilon>0$,
$$
\begin{aligned}
&\E\Big[\tilde{\alpha}^H_\epsilon\big([\delta, 1]
\times [0,1]^{p-1}\big)^m\Big]\\
&=\int_{(\R^d)^m}dx_1\cdots dx_m\bigg[\int_{[\delta, 1]^m}\E\prod_{k=1}^m
p_\epsilon\big(W^H(s_k)-x_k\big)ds_1\cdots ds_m\bigg]\\
&\hskip.4in\times\bigg[\int_{[0, 1]^m}\E\prod_{k=1}^m
p_\epsilon\big(W^H(s_k)-x_k\big)ds_1\cdots ds_m\bigg]^{p-1}\\
&\le\bigg\{\int_{(\R^d)^m}dx_1\cdots dx_m\bigg[\int_{[\delta, 1]^m}\E\prod_{k=1}^m
p_\epsilon\big(W^H(s_k)-x_k\big)ds_1\cdots ds_m\bigg]^p\bigg\}^{1/p}\\
&\times\bigg\{\int_{(\R^d)^m}dx_1\cdots dx_m\bigg[\int_{[0, 1]^m}\E\prod_{k=1}^m
p_\epsilon\big(W^H(s_k)-x_k\big)ds_1\cdots ds_m\bigg]^p\bigg\}^{(p-1)/ p}\\
&=\bigg\{\E\Big[\tilde{\alpha}^H_\epsilon\big([\delta, 1]^p\big)^m\Big]
\bigg\}^{1/p}\bigg\{\E\Big[\tilde{\alpha}^H_\epsilon\big([0, 1]^p\big)^m\Big]
\bigg\}^{(p-1)/p}\, .
\end{aligned}
$$
Letting $\epsilon\to 0^+$ yields
$$
\E\Big[\tilde{\alpha}^H\big([\delta, 1]
\times [0,1]^{p-1}\big)^m\Big]
\le\bigg\{\E\Big[\tilde{\alpha}^H\big([\delta, 1]^p\big)^m\Big]\bigg\}^{1/ p}
\bigg\{\E\Big[\tilde{\alpha}^H\big([0, 1]^p\big)^m\Big]\bigg\}^{(p-1)/p}\, .
$$
Similarly,
$$
\E\Big[\tilde{\alpha}^H\big([0, \delta]
\times [0,1]^{p-1}\big)^m\Big]
\le\bigg\{\E\Big[\tilde{\alpha}^H\big([0, \delta]^p\big)^m\Big]\bigg\}^{1/p}
\bigg\{\E\Big[\tilde{\alpha}^H\big([0, 1]^p\big)^m\Big]\bigg\}^{(p-1)/ p}\, .
$$
So we have
$$
\bigg\{\E\Big[\tilde{\alpha}^H\big([0,1]^p\big)^m\Big]\bigg\}^{1/ mp}\le
\bigg\{\E\Big[\tilde{\alpha}^H\big([\delta, 1]^p\big)^m\Big]\bigg\}^{1/ mp}
+\bigg\{\E\Big[\tilde{\alpha}^H\big([0, \delta]^p\big)^m\Big]\bigg\}^{1/ mp} \, .
$$

By scaling,
$$
\E\Big[\tilde{\alpha}^H\big([0, \delta]^p\big)^m\Big]
=\delta^{(p-Hd(p-1))m}\E\Big[\tilde{\alpha}^H\big([0, 1]^p\big)^m\Big]\, .
$$
Thus
$$
\E\Big[\tilde{\alpha}^H\big([\delta, 1]^p\big)^m\Big]
\ge \Big[1-\delta^{1-Hd(p-1)/p}\Big]^{mp}
\E\Big[\tilde{\alpha}^H\big([0, 1]^p\big)^m\Big]\,.
$$
Therefore, (\ref{5.2.5}) follows from (\ref{eq:limTh4}).

To bound the limit in (\ref{eq:limTh3}) from below, we claim that
\begin{align}\label{5.2.6}
&\lim_{m\to\infty}\frac{1}{m}\log\left\{ (m!)^{-Hd(p-1)}
\E\alpha^{H}\big([0,1]^{p}\big)^{m}\right\}\\
&\ge
p\log\bigg\{(1-Hd/p^*)^{-(1-Hd/p^*)}(p^*)^{d\over 2p^*}(2\pi)^{-{d\over 2p^*}}
\int_0^\infty \big(1+t^{2H}\big)^{-d/2}e^{-t}dt\bigg\}.\nonumber
\end{align}
{}

Let $\tau_1,\cdots, \tau_p$ be i.i.d. exponential times independent
of $B_1^H(t),\cdots, B_p^H(t)$. Given $\epsilon>0$
$$
\begin{aligned}
\E\Big[\alpha_\epsilon^H \big([0,\tau_1]\times\cdots
\times [0,\tau_p]\big)^m\Big]
=\int_{(\R^d)^m}dx_1\cdots dx_m \, Q_\epsilon^p (x_1,\cdots, x_m)
\end{aligned}
$$
where
$$
Q_\epsilon(x_1,\cdots, x_m)
=\int_0^\infty e^{-t}\bigg[\int_{[0,t]^m}ds_1\cdots ds_m
\E\prod_{k=1}^mp_\epsilon \big(B^H(s_k)-x_k\big)\bigg]dt\,.
$$

Let $f(x_1,\cdots, x_m)$ be a rapidly decreasing function on $(\R^d)^m$
such that
$$
\int_{(\R^d)^m}\vert f(x_1,\cdots, x_m)\vert^{p^*}dx_1\cdots dx_m=1\,.
$$
By H\"older inequality,
$$
\begin{aligned}
&\Big\{\E\Big[\alpha_\epsilon^H \big([0,\tau_1]\times\cdots
\times [0,\tau_p]\big)\Big]^m\Big\}^{1/p}\\
&\ge \int_{(\R^d)^m}dx_1\cdots dx_m f(x_1,\cdots, x_m)
Q_\epsilon(x_1,\cdots, x_m)\nonumber\\
&=\int_0^\infty e^{-t}\int_{[0,t]^m}
\bigg[\int_{(\R^d)^m}dx_1\cdots dx_m f(x_1,\cdots, x_m) H_{{\bf s}, \epsilon} (x_1,\cdots, x_m)\bigg] \, ds_1\cdots ds_m \, dt, \nonumber
\end{aligned}
$$
where
$$
H_{{\bf s}, \epsilon} (x_1,\cdots, x_m)
=\E\prod_{k=1}^mp_\epsilon \big(B^H(s_k)-x_k\big)\hskip.2in
{\bf s}=(s_1,\cdots, s_m).
$$
Consider the Fourier transform
$$
\widehat{f}(\lambda_1,\cdots,\lambda_m)
=\int_{(\R^d)^m}dx_1\cdots dx_m f(x_1,\cdots, x_m)
\exp\bigg\{i\sum_{k=1}^m\lambda_k\cdot x_k\bigg\}\,.
$$
It is easy to see that
$$
\widehat{H}_{{\bf s}, \epsilon}(\lambda_1,\cdots,\lambda_m)
=\exp\bigg\{-{\epsilon\over 2}\sum_{k=1}^m\vert\lambda_k\vert^2
-{1\over 2}\Var\Big(\sum_{k=1}\lambda_k\cdot B^H(s_k)\Big)\bigg\}\,.
$$
By Parseval identity,
\begin{align*}
\int_{(\R^d)^m}dx_1\cdots dx_m & f(x_1,\cdots, x_m)
H_{{\bf s}, \epsilon} (x_1,\cdots, x_m)\\
&={1\over (2\pi)^{md}}\int_{(\R^d)^m}d\lambda_1\cdots d\lambda_m
\widehat{f}(\lambda_1,\cdots,\lambda_m)\\
&  \quad  \times\exp\bigg\{-{\epsilon\over 2}\sum_{k=1}^m\vert\lambda_k\vert^2
-{1\over 2}\Var\Big(\sum_{k=1}^m\lambda_k\cdot B^H(s_k)\Big)\bigg\}\,.
\end{align*}
Thus,
\begin{align*}
\Big\{\E\Big[\alpha_\epsilon^H & \big([0,\tau_1]\times\cdots
\times [0,\tau_p]\big)\Big]^m\Big\}^{1/p}\\
&\ge {1\over (2\pi)^{md}}
\int_0^\infty e^{-t}dt\int_{[0,t]^m}ds_1\cdots ds_m\bigg[\int_{(\R^d)^m}
d\lambda_1\cdots d\lambda_m\\
&  \quad  \times\widehat{f}(\lambda_1,\cdots,\lambda_m)\exp\bigg\{
-{\epsilon\over 2}\sum_{k=1}^m\vert\lambda_k\vert^2-{1\over 2}\Var\Big(\sum_{k=1}^m\lambda_k\cdot B^H(s_k)\Big)\bigg\}\bigg]\,.
\end{align*}
We now let $\epsilon\to 0^+$ on the both hand sides. Noticing that the left
hand side falls into an obvious similarity to (\ref{5.1.1}),
\begin{align}\label{4-4}
\Big\{\E\Big[\alpha^H & \big([0,\tau_1]\times\cdots
\times [0,\tau_p]\big)\Big]^m\Big\}^{1/p}\\
&\ge {1\over (2\pi)^{md}}
\int_0^\infty e^{-t}dt\int_{[0,t]^m}ds_1\cdots ds_m\bigg[\int_{(\R^d)^m}
d\lambda_1\cdots d\lambda_m\nonumber\\
&  \quad \times\widehat{f}(\lambda_1,\cdots,\lambda_m)\exp\bigg\{
-{1\over 2}\Var\Big(\sum_{k=1}^m\lambda_k\cdot B^H(s_k)\Big)\bigg\}\bigg]\,.
\nonumber
\end{align}

We now specify the function $f(x_1,\cdots, x_m)$ as
$$
f(x_1,\cdots, x_m)=C^m
\prod_{k=1}^mp_1(x_k)
$$
where
$$
C=(p^*)^{d\over 2p^*}(2\pi)^{d(p^*-1)\over 2p^*}.
$$
We have
$$
\begin{aligned}
&\int_{(\R^d)^m}d\lambda_1\cdots d\lambda_m
\widehat{f}(\lambda_1,\cdots,\lambda_m)\exp\bigg\{
-{1\over 2}\Var\Big(\sum_{k=1}^m\lambda_k\cdot B^H(s_k)\Big)\bigg\}\\
&=C^m
\bigg[\int_{\R^m}d\gamma_1\cdots d\gamma_m
\exp\bigg\{-{1\over 2}\sum_{k=1}^m\gamma_k^2
-{1\over 2}\Var\Big(\sum_{k=1}^m\gamma_k B_0^H(s_k)\Big)\bigg\}\bigg]^d
\end{aligned}
$$
where $B_0^H(t)$ is an 1-dimensional fractional Brownian motion.

Let $\xi_1,\cdots \xi_m$ be i.i.d. standard normal random variable
independent of $B_0^H(t)$. Write
$$
\eta_k=\xi_k+B_0^H(s_k)\hskip.2in k=1,\cdots, m.
$$
We have
$$
{1\over 2}\sum_{k=1}^m\gamma_k^2
+{1\over 2}\Var\Big(\sum_{k=1}^m\gamma_k B_0^H(s_k)\Big)
={1\over 2}\Var\Big(\sum_{k=1}^m\gamma_k \eta_k\Big)
$$
And thus by Gaussian integration,
$$
\begin{aligned}
&\int_{\R^m}d\gamma_1\cdots d\gamma_m
\exp\bigg\{-{\sigma^2\over 2}\sum_{k=1}^m\gamma_k^2
-{1\over 2}\Var\Big(\sum_{k=1}^m\gamma_k B_0^H(s_k)\Big)\bigg\}\\
&=(2\pi)^{m/2}\det\Big\{\Cov\big(\eta_1,\cdots, \eta_m\big)\Big\}^{-1/2}.
\end{aligned}
$$
with convention that $s_0=0$.

Write $s_0=0$ and assume $s_1<\cdots <s_m$. By Lemma \ref{lem:Cov},
$$
\begin{aligned}
&\det\Big\{\Cov\big(\eta_1,\cdots, \eta_m\big)\Big\}
=\Var(\eta_1)\prod_{k=2}^m\Var\Big(\eta_k|\eta_1,\cdots,\eta_{k-1}\Big)\\
&=\Big\{1+\Var\big(B_0(s_1)\big)\Big\}\prod_{k=2}^m
\Big\{1+\Var\Big(B_0^H(s_k)|B_0^H(s_1),\cdots,
B_0^H(s_{k-1})\Big)\Big\}\\
&\le\prod_{k=1}^m\Big\{1+(s_k-s_{k-1})^{2H}\Big\}
\end{aligned}
$$
where the last step follow from the computation
$$
\begin{aligned}
&\Var\Big(B_0^H(s_k)|B_0^H(s_1),\cdots,
B_0^H(s_{k-1})\Big)\\
&=\Var\Big(B_0^H(s_k)-B_0^H(s_{k-1})|B_0^H(s_1),\cdots,
B_0^H(s_{k-1})\Big)\\
&\le \Var\Big(B_0^H(s_k)-B_0^H(s_{k-1})\Big)
=(s_k-s_{k-1})^{2H}.
\end{aligned}
$$
Summarizing our argument since (\ref{4-4}), we obtain
$$
\begin{aligned}
&\Big\{\E\Big[\alpha^H \big([0,\tau_1]\times\cdots
\times [0,\tau_p]\big)\Big]^m\Big\}^{1/p}\\
&\ge m! \big({C(2\pi)^{-d/2}}\big)^m
\int_0^\infty e^{-t}dt\int_{[0,t]_<^m}ds_1\cdots ds_m
\prod_{k=1}^m\Big\{1+(s_k-s_{k-1})^{2H}\Big\}^{-d/2}\\
&= m!\big({C(2\pi)^{-d/2}}\big)^m
\bigg[\int_0^\infty \big(1+t^{2H}\big)^{-d/2}e^{-t}dt\bigg]^m.
\end{aligned}
$$
Equivalently,
\begin{align}\label{4-5}
&\E\Big[\alpha^H \big([0,\tau_1]\times\cdots
\times [0,\tau_p]\big)\Big]^m \ge  (m!)^p \big({C (2\pi)^{-d/2}}\big)^{mp}
\bigg[\int_0^\infty \big(1+t^{2H}\big)^{-d/2}e^{-t}dt
\bigg]^{pm}.
\end{align}
On the other hand, with obvious similarity to (\ref{5.1.7'})
$$
\begin{aligned}
&\E\Big[\alpha^H \big([0,\tau_1]\times\cdots
\times [0,\tau_p]\big)\Big]^m
\le\E\Big[\alpha^H \big([0,1]^p\big)\Big]^m
\Big\{\Gamma\big(1+m(1-Hd/p^*\big)\Big\}^p.
\end{aligned}
$$
Hence, (\ref{5.2.6}) follows from (\ref{4-5}) and Stirling
formula.

By (\ref{5.1.8}) and (\ref{5.2.6}), the limit given in
(\ref{eq:limTh3}) is finite. By Lemma \ref{lem:KM}, the
large deviation given in (\ref{eq:lim-aB}) holds with
$$
\begin{aligned}
K(H, d, p)&=Hd(p-1)\exp\Big\{-{C(H,d,p)-d(p-1)\log c_H\over Hd(p-1)}\Big\}
%\\&=c_H^{1/H}Hd(p-1)\exp\Big\{-{C(H,d,p)\over Hd(p-1)}\Big\}
=c_H^{1/H}\tilde{K}(H, d, p).
\end{aligned}
$$
\hfill\qed

\section{The law of the iterated logarithm} \label{s:LIL}

We will prove Theorem \ref{1-17} in this section. Due to the similarity of arguments, we will only establish (\ref{1-21}). By the self-similarity property \eqref{eq:aW-ss}, the large deviation limit of Theorem \ref{th:4} can be rewritten as
\begin{align}\label{6.1}
&\lim_{t\to\infty}(\log\log t)^{-1}\log\P
\Big\{\tilde{\alpha}^H\big([0,t]^p\big)\ge\lambda t^{p-Hd(p-1)}
(\log\log t)^{Hd(p-1)}\Big\}\nonumber\\
&=-\tilde{K}(H,d, p)\lambda^{p^*/Hdp}\hskip.2in (\lambda>0).
\end{align}
Therefore, the upper bound
$$
\limsup_{t\to\infty} t^{Hd(p-1)-p}(\log\log t)^{-Hd(p-1)}
\widetilde{\alpha}^{H}\big([0,t]^{p}\big)
\le\widetilde{K}(H, d, p)^{-Hd(p-1)}
\hskip.2in a.s
$$
is a consequence of the standard argument using Borel-Cantelli lemma.

To show the lower bound, we proceed in several steps. First
let $N > 1$ be a large but fixed number and write $t_n=N^n$
($n=1,2,\cdots$).  Define the $d$-dimensional
process
$$
Q_n^H(t)=\int_0^{t_n}(t+u)^{H-1/2}dB(u)\hskip.2in t\ge 0,
$$
where $B(u)$ is a standard $d$-dimensional Brownian motion. Recall that $\mathbb{H}[0,T]$ denotes the RKHS of $\{W^H(t)\}_{t\in [0,T]}$. Combining Propositions \ref{pro:rem} and  \ref{prop:Z_a} we can deduce that $\{Q_n^H(t)\}_{t\in [0,T]}$ is not in $\mathbb{H}[0,T]$, $T>0$. For that reason, similarly as in Proposition \ref{prop:Z_a}, we define the following modifications of $Q_n^H(t)$.
When $H\in (0, 1/2)$, put
$$
\begin{aligned}
G_n^{H}(t)=\begin{cases}
A_n t, & 0\le t\le t_n\\
Q_n^{H}(t), & t>t_n,\end{cases}
\end{aligned}
$$
where $A_n=t_n^{-1}Q_n^{H}(t_n)$. 
When $H\in(\frac{1}{2},1)$, put
$$
\begin{aligned}
G_n^{H}(t)=\begin{cases}
B_{1,n} t^{2} + B_{2,n} t^{3}, & 0\le t\le t_n\\
Q_n^{H}(t), & t>t_n, \end{cases}
\end{aligned}
$$
where $B_{1,n}=3t_n^{-2}Q_n^{H}(t_n)-t_n^{-1}\dot{Q}_n^{H}(t_n)$ and $B_{2,n}=  -2t_n^{-3}Q_n^H(t_n)+t_n^{-2}\dot{Q}_n^{H}(t_n)$.

%:here2

\begin{lemma}\label{6.2}
For every $n\ge 1$, $\P \left(\{G_n^H(t)\}_{t\in [0,t_{n+1}]} \subset \mathbb{H}[0,t_{n+1}] \right)=1$. Furthermore,
\begin{align}\label{6.3}
\sup_n \E\| G_n^H\|_{\mathbb{H}[0,t_{n+1}]}^2<\infty.
\end{align}
\end{lemma}

\proof
Obviously, it suffices to consider the case $d=1$. The first part of the lemma follows by the same argument as in Proposition \ref{prop:Z_a}. For the second part we use Lemma \ref{lem:bound} with $a=t_n$ and $T=t_{n+1}$. A constant $C>0$ below will depend only on $H$ but it will be allowed to be different at different places. 

First consider $H\in (0, 1/2)$, so that $m=\lceil H+1/2\rceil =1$. In this case we get
\begin{align*}
\|\dot{G}_n^{H}\|_{L_{\infty}[0,t_n]} = |A_n|
\end{align*}
and 
\begin{align*}
\int_{t_n}^{t_{n+1}} & \left|\int_{t_n}^t (t-s)^{-H-1/2}\dot{G}_n^{H}(s) \, ds \right|^2 dt \\
&= C\int_{t_n}^{t_{n+1}}  \left|\int_0^{t_n}  \left( \int_{t_n}^t (t-s)^{-H-1/2} (s+u)^{H-3/2} \, ds \right) dB(u) \right|^2 dt.
\end{align*}
Therefore,
\begin{equation}\label{eq:b1}
\E \|\dot{G}_n^{H}\|_{L_{\infty}[0,t_n]}^2 = t_n^{-2} \E Q_n^{H}(t_n)^2 = C t_n^{2H-2}, 
\end{equation}
and
\begin{align}
\E \int_{t_n}^{t_{n+1}} & \left|\int_{t_n}^t (t-s)^{-H-1/2}\dot{G}_n^{H}(s) \, ds \right|^2 dt \nonumber \\
&= C\int_{t_n}^{t_{n+1}}  \int_0^{t_n}  \left( \int_{t_n}^t (t-s)^{-H-1/2} (s+u)^{H-3/2} \, ds \right)^2 du \, dt \nonumber \\
& = C\int_{t_n}^{t_{n+1}}  \int_0^{t_n}  \frac{(t-t_n)^{1-2H}(u+t_n)^{2H-1}}{(t+u)^2}  du \, dt \nonumber \\
&\le C \int_{t_n}^{t_{n+1}}  \left( \frac{t}{t_n} -1\right)^{1-2H} \frac{t_n}{t^2}\, dt \le C \left( \frac{t_{n+1}}{t_n}-1\right)^{1-2H} \int_{t_n}^{t_{n+1}}\frac{1}{t} \, dt \nonumber \\
& \le  C \left( \frac{t_{n+1}}{t_n}\right)^{1-2H}\log \frac{t_{n+1}}{t_{n}}\, . \label{eq:b2}
\end{align}
Using bounds \eqref{eq:b1}-\eqref{eq:b2} with Lemma \ref{lem:bound} we get 
\begin{align*}
\E\| G_n^H\|_{\mathbb{H}[0,t_{n+1}]}^2 &\le C (t_{n+1}^{2-2H}-t_n^{2-2H}) t_n^{2H-2} + C (t_{n+1}/t_n)^{1-2H}\log(t_{n+1}/t_n)\\
& \le C( N^{2-2H} + N^{1-2H}\log N)\, ,
\end{align*}
which proves \eqref{6.3} in the case $H \in (0,1/2)$. The proof in the case $H \in (1/2,1)$ follows the same line of computations, thus is omitted.
\hfill \qed

\medskip

For a simplicity of notation, from now on write $\mathbb{H}_n$ for $\mathbb{H}[0,t_{n+1}]$.
Define the sigma field
$$
{\cal F}_t=\sigma\Big\{\big(B_1(s),\cdots, B_p(s)\big);
\hskip.1in s\le t\Big\}.
$$
To complete the proof of Theorem \ref{1-17}, i.e., to establish the lower bound in \eqref{1-21}, it is enough to show that for any $\lambda<\widetilde{K}(H, d, p)^{-Hd(p-1)}$ there is an $N$, sufficiently large, such that
\begin{align}\label{6.4}
\sum_n\P\Big\{\tilde{\alpha}^H\big([2t_n, t_{n+1}]^p\big)\ge\lambda
t_{n+1}^{p-Hd(p-1)}(\log\log t_{n+1})^{Hd(p-1)}\Big\vert {\cal F}_{t_n}\Big\}
=\infty \hskip.2in a.s.
\end{align}
Indeed, by  \cite[Corollary 5.29 p. 96]{Breiman}, (\ref{6.4}) implies that
$$
\limsup_{n\to\infty} t_{n+1}^{Hd(p-1)-p}(\log\log t_{n+1})^{-Hd(p-1)}
\tilde{\alpha}^H \big([2t_n, t_{n+1}]^p\big)
\ge \lambda\hskip.2in a.s.
$$
which leads to
$$
\limsup_{t\to\infty} t^{Hd(p-1)-p}(\log\log t)^{-Hd(p-1)}
\tilde{\alpha}^H \big([0, t]^p\big)
\ge \lambda\hskip.2in a.s.
$$
Letting $\lambda\to\widetilde{K}(H, d, p)^{-Hd(p-1)}$ on the right hand
side yields the lower bound as claimed.

Let now $\epsilon>0$ be fixed and write
$$
\begin{aligned}
&\tilde{\alpha}_\epsilon^H\big([2t_n, t_{n+1}]^p\big)
=\int_{[2t_n, t_{n+1}]^p}ds_1\cdots ds_p \, g_\epsilon \big(W_1^H(s_1),\cdots,
W_p^H(s_p)\big)\\
&=\int_{[t_n, t_{n+1}-t_n]^p}ds_1\cdots ds_p \, g_\epsilon \big(W_1^H(t_n+s_1),\cdots,
W_p^H(t_n+s_p)\big)\\
&=\int_{[t_n, t_{n+1}-t_n]^p}ds_1\cdots ds_p \, g_\epsilon
\big(Y_1^H(s_1)+Z_1^H(s_1),\cdots, Y_p^H(s_p)+Z_p^H(s_p)\big),
\end{aligned}
$$
where $g_\epsilon(x_1,\cdots, x_p)$ is given in (\ref{eq:g_e}) and
for $j=1,\cdots, p$,
$$
Y_j^H(t)=\int_{t_n}^{t_n+t}(t_n+t-s)^{H-1/2}dB_j(s),\hskip.2in
Z_j^H(t)=\int_0^{t_n}(t_n+t-s)^{H-1/2}dB_j(s)\,.
$$
Consider a symmetric set
$A\subset\otimes_{j=1}^p C\Big\{[0, t_{n+1}], \R^d\Big\}$
defined by
$$
\begin{aligned}
&A=\Bigg\{(f_1,\cdots, f_p)\in\otimes_{j=1}^p C\Big\{[0, t_{n+1}], \R^d\Big\};
\\
&\int_{[t_n, t_{n+1}-t_n]^p}ds_1\cdots ds_p \,
g_\epsilon \big(f_1(s_1),\cdots, f_p(s_p)\big)\ge\lambda
t_{n+1}^{p-Hd(p-1)}(\log\log t_{n+1})^{Hd(p-1)}\Bigg\}\,.
\end{aligned}
$$

For any $(f_1,\cdots, f_p)\in \otimes_{j=1}^p \mathbb{H}_n$, applying
Lemma \ref{lem:A-CM}-(ii) to the indicator of $A$ leads to
$$
\P\Big\{\big(W_1^H+f_1,\cdots, W_p^H+f_p)\in A\Big\}
\ge \exp\Big\{-{1\over 2}\sum_{j=1}^p \| f\|_{H_n}^2\Big\}
\P\Big\{\big(W_1^H,\cdots, W_p^H)\in A\Big\},
$$
if $f_1,\cdots, f_p\in \mathbb{H}_n$.

Notice that
$$
\Big\{Z^H(t);\hskip.1in t_n\le t\le t_{n+1}\Big\}\buildrel d\over
=\Big\{Q_n^H(t);\hskip.1in t_n\le t\le t_{n+1}\Big\}
=\Big\{G_n^H(t);\hskip.1in t_n\le t\le t_{n+1}\Big\}
$$
$$
\Big\{Y^H(t);\hskip.1in t_n\le t\le t_{n+1}\Big\}\buildrel d\over
=\Big\{W^H(t);\hskip.1in t_n\le t\le t_{n+1}\Big\}
$$
and $Y^H(t)$ and $Z^H(t)$ are independent. By Lemma \ref{6.2},
\begin{align*}
\P\Big\{ \big(Y_1^H+Z^H_1,\cdots, & Y_p^H+Z^H_p) 
\in A\Big\vert {\cal F}_{t_n}\Big\} \\
&\ge \exp\Big\{-{1\over 2}\sum_{j=1}^p \| G_{n, j}^H \|_{\mathbb{H}_n}^2
\Big\}\P\Big\{\big(W_1^H,\cdots, W_p^H)\in A\Big\},
\end{align*}
or
$$
\begin{aligned}
&\P\Big\{\tilde{\alpha}_\epsilon^H\big([2t_n, t_{n+1}]^p\big)\ge\lambda
t_{n+1}^{p-Hd(p-1)}(\log\log t_{n+1})^{Hd(p-1)}\Big\vert {\cal F}_{t_n}\Big\}\\
&\ge \exp\Big\{-{1\over 2}\sum_{j=1}^p \| G_{n, j}^H\|_{\mathbb{H}_n}^2
\Big\}\P\Big\{\tilde{\alpha}_\epsilon^H\big([t_n, t_{n+1}-t_n]^p\big)\ge\lambda
t_{n+1}^{p-Hd(p-1)}(\log\log t_{n+1})^{Hd(p-1)}\Big\}.
\end{aligned}
$$
Letting $\epsilon\to 0^+$ on the both sides yields
$$
\begin{aligned}
&\P\Big\{\tilde{\alpha}^H\big([2t_n, t_{n+1}]^p\big)\ge\lambda
t_{n+1}^{p-Hd(p-1)}(\log\log t_{n+1})^{Hd(p-1)}\Big\vert {\cal F}_{t_n}\Big\}\\
&\ge \exp\Big\{-{1\over 2}\sum_{j=1}^p \| G_{n, j}^H\|_{\mathbb{H}_n}^2
\Big\}\P\Big\{\tilde{\alpha}^H\big([t_n, t_{n+1}-t_n]^p\big)\ge\lambda
t_{n+1}^{p-Hd(p-1)}(\log\log t_{n+1})^{Hd(p-1)}\Big\}.
\end{aligned}
$$

By (\ref{6.1}) and by an argument similar to the one used for
(\ref{5.2.5}), for $\lambda<\widetilde{K}(H, d, p)^{-Hd(p-1)}$
and any small $\delta>0$,
one can take $N$ sufficiently large so that, for large $n$,
\begin{align*}
\P\Big\{\tilde{\alpha}^H\big([t_n, t_{n+1}-t_n]^p \big)  & \ge \lambda
t_{n+1}^{p-Hd(p-1)}(\log\log t_{n+1})^{Hd(p-1)}\Big\}\\
&\ge
\exp\big\{-(1-\delta)\log\log t_{n+1}\big\}=(n\log N)^{-1+\delta}.
\end{align*}
To establish (\ref{6.4}), therefore, it suffices to show that
 for any $\epsilon,\delta>0$,
\begin{align}\label{6.5}
\sum_n{1\over n^{1-\delta}}\1\Big\{\sum_{j=1}^p \| G_{n,j}^H
\|_{\mathbb{H}_n}^2
\le\epsilon\log\log t_{n+1}\Big\}=\infty \hskip.2in a.s.
\end{align}

Indeed, by Lemma \ref{6.2} $G_n^H$ can be viewed as a Gaussian
sequence taking values in $H_n$. By the Gaussian tail
estimate, see \cite{LT}, p.59,
there is $u=u(\epsilon)>0$ such that 
$$
\P\Big\{\sum_{j=1}^p \| G_{n, j}^H\|_{\mathbb{H}_n}^2
\ge\epsilon\log\log t_{n+1}\Big\}\le {1\over n^u}
$$
for large $n$. Then for $0<\delta<u$, we obtain (\ref{6.5}), 
which yields (\ref{6.4}). The proof is complete.
 \hfill\qed

\section{Local times of Gaussian fields}
\label{construction}

We begin with mentioning the work of Geman, Horowitz and Rosen
(\cite{GHR}) on the condition for the existence
and continuity of the local times of the Gaussian fields, see also recent work of Wu and Xiao \cite{WX}.
Let $X({\bf t})$ (${\bf t}\in(\R^+)^p$) be a mean zero Gaussian
field taking values in $\R^d$ such that there is a $\gamma>0$ such that
for any $t>0$ and $m \in \N$, 
%$m=1,2,\cdots$,
\begin{align}\label{7.1}
\int_{([0, t]^p)^m}  d{\bf s}_1\cdots d{\bf s}_m &
 \int_{(\R^d)^m}  d\lambda_1\cdots d\lambda_m\\
&  \quad \times\Big(\prod_{k=1}^m\vert\lambda_k\vert^\gamma\Big)
\exp\bigg\{-{1\over 2}
\Var\Big(\sum_{k=1}^m
\lambda_k\cdot X({\bf s}_k)\Big)\bigg\}<\infty.\nonumber
\end{align}
 Geman, Horowitz and Rosen
(Theorem (2.8) in \cite{GHR}) proved that the occupation time
$$
\mu_{\bf t}(B)=\int_{[{\bf 0}, {\bf t}]}\1_{\{X({\bf s})\in B\}} \, d{\bf s}\hskip.2in
B\subset\R^d
$$
is absolutely continuous with respect to the Lebesgue measure on $\R^d$.
Further, the correspondent density function formally written as
$$
\alpha \big([{\bf 0},{\bf t}], x\big)
=\int_{[{\bf 0}, {\bf t}]}\delta_x\big(X({\bf s})\big) \, d{\bf s}
$$
is jointly continuous in $({\bf t}, x)$. For fixed $x$,
the distribution function $\alpha \big([{\bf 0},{\bf t}], x\big)$
(${\bf t}\in (\R^+)^p$) generates a (random) measure
$\alpha (A, x)$ ($A\subset (\R^+)^p$) on $(\R^+)^p$ which is called
the local time of $X({\bf t})$.

In this paper, the result of Geman, Horowitz and Rosen is applied
to the following four Gaussian fields:

1. The $d$-dimensional fractional Brownian motion $X_1(t)=B^H(t)$.

2. The $d$-dimensional Riemann-Liouville process $X_2(t)=W^H(t)$.

3. The $d(p-1)$-dimension Gaussian field
$$
X_3(t_{1},\cdots,t_{p})
=\Big(B_{1}^{H}(t_{1})-B_{2}^{H}(t_{2}),
\cdots,B_{p-1}^{H}(t_{p-1})-B_{p}^{H}(t_{p})\Big).
$$
4. The $d(p-1)$-dimension Gaussian field
$$
X_4(t_{1},\cdots,t_{p})
=\Big(W_{1}^{H}(t_{1})-W_{2}^{H}(t_{2}),
\cdots,W_{p-1}^{H}(t_{p-1})-W_{p}^{H}(t_{p})\Big).
$$

\begin{theorem}\label{7.2}
Under $Hd<1$, $X_1(t)$ and $X_2(t)$
satisfy the condition (\ref{7.1}); under $Hd<p^*$, $X_3({\bf t})$ and
$X_4({\bf t})$
satisfy the condition (\ref{7.1}). Consequently, $X_1$, $X_2$,
$X_3$ and $X_4$ have continuous (jointly in time and space variables)
local times.
\end{theorem}

\proof Due to similarity we only verify (\ref{7.1}) for $X_3$, which becomes
\begin{align}\label{6-1}
 & \int_{([0,t]^{p})^{m}}d{\bf s}_{1}\cdots d{\bf s}_{m}
\int_{(\R^{d(p-1)})^{m}}d\tilde{\lambda}_{1}\cdots d\tilde{\lambda}_{m}\exp\bigg\{-\frac{1}{2}\Var\Big(\sum_{k=1}^{m}
\tilde{\lambda}_{k}\cdot X({\bf s}_{k})\Big)\bigg\}
\prod_{k=1}^{m}\vert\tilde{\lambda}_{k}\vert^{\gamma}<\infty, 
 \end{align}
where we use the notation
$$
{\bf s}_k=(s_{k,1},\cdots, s_{k,p})\hskip.1in\hbox{and}\hskip.1in
\tilde{\lambda}_{k}=(\lambda_{k,1},\cdots,\lambda_{k,p-1})\,.
$$

Notice that
\[
\Var\Big(\sum_{k=1}^{m}\tilde{\lambda}_{k}\cdot X({\bf s}_{k})\Big)
=\sum_{j=1}^{p}\Var\Big(\sum_{k=1}^{m}(\lambda_{k,j}-\lambda_{k,j-1})\cdot B^{H}(s_{k,j})\Big)
\]
with the convention
$\lambda_{k,0}=\lambda_{k,p}=0$.
By suitable substitution and using the bound
\[
\vert\tilde{\lambda}_{k}\vert\le
C\prod_{j=1}^{p}\max\{1,\vert\lambda_{k,j}-\lambda_{k,j-1}\vert\},
\]
 we have
 \[
\begin{aligned}
& \int_{(\R^{d(p-1)})^{m}}d\tilde{\lambda}_{1}\cdots d\tilde{\lambda}_{m}
\exp\bigg\{-\frac{1}{2}\Var\Big(\sum_{k=1}^{m}\tilde{\lambda}_{k}\cdot X({\bf s}_{k})\Big)\bigg\}\prod_{k=1}^{m}\vert\tilde{\lambda}_{k}\vert^{\gamma}\\
 & \le C\int_{(\R^{md})^{p-1}}d\bar{\lambda}_{1}\cdots
d\bar{\lambda}_{p-1}\prod_{j=1}^{p}H_{j}(\bar{\lambda}_{j})
 \end{aligned}
\]
 where
 \[
H_{j}(\bar{\lambda}_j)=\Big(\prod_{k=1}^{m}
\max\{1,\vert\lambda_{k, j}\vert^{\gamma}\}\Big)\exp\bigg\{-\frac{1}{2}\Var\Big(\sum_{k=1}^{m}\lambda_{k, j}\cdot B^{H}(s_{k,j})\Big)\bigg\}
\]
for $\bar{\lambda}_j=(\lambda_{1, j},\cdots,\lambda_{m, j})$ ($1\le j\le p-1$)
and $\bar{\lambda}_{p}=-(\bar{\lambda}_{1}+\cdots+\bar{\lambda}_{p-1})$.

Write
$$
\prod_{j=1}^pH_j(\bar{\lambda}_j)=\prod_{j=1}^p\prod_{1\le k\neq j \le p}H_k(\bar{\lambda}_k)^{1/(p-1)}.
$$
By H\"older inequality
\begin{align*}
\int_{(\R^{md})^{p-1}}d\bar{\lambda}_1\cdots
d\bar{\lambda}_{p-1}\prod_{j=1}^pH_j(\bar{\lambda}_j)
\le\prod_{j=1}^p\bigg\{\int_{(\R^{md})^{p-1}}d\bar{\lambda}_1\cdots
d\bar{\lambda}_{p-1}\prod_{1\le k\neq j \le p}H_k(\bar{\lambda}_k)^{p^*}\bigg\}^{1/p}.
\end{align*}

When $j=p$,
$$
\begin{aligned}
&\int_{(\R^{md})^{p-1}}d\bar{\lambda}_1\cdots
d\bar{\lambda}_{p-1}\prod_{1\le k< p}H_k(\bar{\lambda}_k)^{p^*}
=\prod_{k=1}^{p-1}\int_{\R^{md}}H_k(\bar{\lambda})^{p^*}d\bar{\lambda}\,.
\end{aligned}
$$

As for $1\le j\le p-1$, recall that
$\bar{\lambda}_p=-(\bar{\lambda}_1+\cdots +\bar{\lambda}_{p-1})$.
By translation invariance,
$$
\int_{\R^{md}}H_p(\bar{\lambda}_p)^{p^*}d\bar{\lambda}_j=
\int_{\R^{md}}H_p(\bar{\lambda})^{p^*}d\bar{\lambda}.
$$
By Fubini theorem, for fixed $j$,
$$
\begin{aligned}
&\int_{(\R^{md})^{p-1}}d\bar{\lambda}_1\cdots
d\bar{\lambda}_{p-1}\prod_{1\le k\neq j \le p}H_k(\bar{\lambda}_k)^{p^*}
=\prod_{1\le k\neq j \le p}
\int_{\R^{md}}H_k(\bar{\lambda})^{p^*}d\bar{\lambda}.
\end{aligned}
$$
%Summarizing what we have,

Summarize our argument, the left hand of (\ref{6-1}) is bounded by 
\[
\begin{aligned}
%& \int_{([0,t]^{p})^{m}}d{\bf s}_{1}
%\cdots d{\bf s}_{m}\int_{(\R^{d(p-1)})^{m}}d\bar{\lambda}_{1}
%\cdots d\bar{\lambda}_{m}\\
% & \times\exp\bigg\{-\frac{1}{2}\Var\Big(\sum_{k=1}^{m}\bar{\lambda}_{k}\cdot X({\bf s}_{k})\Big)\bigg\}\prod_{k=1}^{m}\vert\bar{\lambda}_{k}\vert^{\gamma}\\
 & C\Bigg\{\int_{[0,t]^{m}}ds_{1}\cdots ds_{m}\bigg[\int_{(\R^{d})^{m}}
d{\lambda}_{1}\cdots d{\lambda}_{m}\Big(\prod_{k=1}^{m}\max\{1,
\vert{\lambda}_{k}\vert^{p^*\gamma}\}\Big)\\
 & \times\exp\bigg\{-\frac{p^*}{2}\Var\Big(\sum_{k=1}^{m}{\lambda}_{k}\cdot B^{H}(s_{k})\Big)\bigg\}\bigg]^{1/p^*}\Bigg\}^{p}.
 \end{aligned}
\]

Hence all we need is to find $\gamma>0$ such
that
\begin{equation}\label{6-2}
  \int_{[0,t]^{m}}ds_{1}\cdots ds_{m}\bigg[\int_{(\R^{d})^{m}}
d{\lambda}_{1}\cdots d{\lambda}_{m}\Big(\prod_{k=1}^{m}\vert{\lambda}_{k}
\vert^{\gamma}\Big)\exp\bigg\{-\frac{p^*}{2}\Var\Big(\sum_{k=1}^{m}
{\lambda}_{k}\cdot B^{H}(s_{k})\Big)\bigg\}\bigg]^{1/p^*}<\infty
 \end{equation}
 for all $m\in \N$.
Further separating variable and substituting variable, the above is
reduced to
\begin{align}\label{6-3}
\int_{[0,t]^{m}}ds_{1}\cdots ds_{m}\bigg[\int_{\R^{m}}d\lambda_{1}\cdots d\lambda_{m}\Big(\prod_{k=1}^{m}\vert\lambda_{k}\vert^{\gamma}\Big) \exp\bigg\{-\frac{1}{2}\Var\Big(\sum_{k=1}^{m}\lambda_{k}B_{0}^{H}(s_{k})\Big)\bigg\}\bigg]^{d/p^*}<\infty .
 \end{align}

By (\ref{3-7}), for any $0=s_0<s_{1}<\cdots<s_{k}$,
\[
\begin{aligned}
& \Var\big(B_{0}^{H}(s_{k})-B_{0}^{H}(s_{k-1})|B_{0}^{H}(s_{1}),\cdots,B_{0}^{H}(s_{k-1})\big)\\
 & \ge\frac{1}{2H}(s_{k}-s_{k-1})^{2H}=\frac{1}{2H}
\Var\big(B_{0}^{H}(s_{k})-B_{0}^{H}(s_{k-1})\big).
 \end{aligned}
\]
 This property is generalized into the notion known as local non-determinism.
By Lemma 2.3 in Berman \cite{Berman}, there is constant $c_{m}>0$
such that for any $\lambda_{1},\cdots,\lambda_{m}\in\R$ and any $s_{1}<\cdots<s_{m}$
\[
\Var\Big(\sum_{k=1}^{m}\lambda_{k}\big(B_{0}^{H}(s_{k})-B_{0}^{H}(s_{k-1})\big)\Big)\ge c_{m}\sum_{k=1}^{m}(s_{k}-s_{k-1})^{2H}\lambda_{k}^{2} \, .
\]
 Consequently, with notation $\lambda_0=0$,
 \[
\begin{aligned}
& \int_{\R^{m}}d\lambda_{1}\cdots d\lambda_{m}\Big(\prod_{k=1}^{m}\vert\lambda_{k}\vert^{\gamma}\Big)\exp\bigg\{-\frac{1}{2}\Var\Big(\sum_{k=1}^{m}\lambda_{k}B_{0}^{H}(s_{k})\Big)\bigg\}\\
 & =\int_{\R^{m}}d\lambda_{1}\cdots d\lambda_{m}\Big(\prod_{k=1}^{m}\vert\lambda_{k}-\lambda_{k-1}\vert^{\gamma}\Big)\\
 & \times\exp\bigg\{-\frac{1}{2}\Var\Big(\sum_{k=1}^{m}\lambda_{k}\big(B_{0}^{H}(s_{k})-B_{0}^{H}(s_{k-1})\big)\Big)\bigg\}\\
 & \le\int_{\R^{m}}d\lambda_{1}\cdots d\lambda_{m}\Big(\prod_{k=1}^{m}\vert\lambda_{k}-\lambda_{k-1}\vert^{\gamma}\Big) \exp\bigg\{-c_{m}\sum_{k=1}^{m}(s_{k}-s_{k-1})^{2H}\lambda_{k}^{2}\bigg\} \, .
 \end{aligned}
\]
 Using triangle inequality (for which we take $\gamma\le 1$)
 \[
\prod_{k=1}^{m}\vert\lambda_{k}-\lambda_{k-1}\vert^{\gamma} \le \prod_{k=1}^{m}\big(\vert\lambda_{k}\vert^\gamma+\vert\lambda_{k-1}\vert^\gamma\big)
=\sum_{j_1,\cdots, j_m}\prod_{k=1}^m\vert\lambda_{k}\vert^{\delta_{j_k}},
\]
where $\delta_{j_k}=0$, $\gamma $ or $2\gamma$. Notice that
$$
\prod_{k=1}^m\vert\lambda_{k}\vert^{\delta_{j_k}}
\le \prod_{k=1}^m(1\vee\vert\lambda_{k}\vert)^{\delta_{j_k}}
\le\prod_{k=1}^m(1\vee\vert\lambda_{k}\vert)^{2\gamma} \,.
$$
Notice the number of the terms in the previous summation is
at most $2^m$. Thus,
$$
\prod_{k=1}^{m}\vert\lambda_{k}-\lambda_{k-1}\vert^{\gamma}
\le 2^m\prod_{k=1}^m(1\vee\vert\lambda_{k}\vert)^{2\gamma}\,.
$$
In this way, the problem is reduced to finding $\gamma>0$ such
that
 \begin{align}
\int_{[0,t]_{<}^{m}} ds_{1}\cdots ds_{m}\,  \bigg[\int_{\R^{m}}d\lambda_{1}\cdots d\lambda_{m}\Big(\prod_{k=1}^{m}\vert\lambda_{k}\vert^{\gamma}\Big) \exp\bigg\{-c_{m}\sum_{k=1}^{m}(s_{k}-s_{k-1})^{2H}\lambda_{k}^{2}\bigg\}\bigg]^{d/p^*}<\infty\,.  \label{6-4}
 \end{align}
Observe that
\[
\begin{aligned}
 \int_{\R^{m}}d  \lambda_{1}\cdots d\lambda_{m} \, & \Big(\prod_{k=1}^{m}\vert\lambda_{k}\vert^{\gamma}\Big)\exp\bigg\{-c_{m}\sum_{k=1}^{m}(s_{k}-s_{k-1})^{2H}\lambda_{k}^{2}\bigg\}  =\prod_{k=1}^{m}\int_{-\infty}^{\infty}\vert\lambda\vert^{\gamma}e^{-c_{m}(s_{k}-s_{k-1})^{2H}\lambda^{2}}d\lambda \\
 &=\bigg\{\int_{-\infty}^{\infty}\vert\lambda\vert^{\gamma}e^{-c_{m}\lambda^{2}}d\lambda\bigg\}^{m}\prod_{k=1}^{m}(s_{k}-s_{k-1})^{-(1+\gamma)H}\,.
 \end{aligned}
\]
Therefore, we need to choose $\gamma>0$ such that
\[
\int_{[0,t]_{<}^{m}}ds_{1}\cdots ds_{m}
\prod_{k=1}^{m}(s_{k}-s_{k-1})^{-(1+\gamma)Hd/p^*}<\infty.
\]
This is always possible because $Hd< p^{\ast}$, so that $(1+\gamma)Hd<p^*$ for some $\gamma>0$. The proof is complete.
\hfill \qed

%\appendix
%%dummy comment inserted by tex2lyx to ensure that this paragraph is not empty
\setcounter{equation}{0} \setcounter{theorem}{0} \global\long\global\long\def\theequation{A\arabic{equation}}
 \global\long\global\long\def\thetheorem{A\arabic{theorem}}

\section{Appendix }

\begin{lemma}\label{constant}
Let $\{B^{H}(t)\}_{t\in\R}$ be a
standard fractional Brownian motion given by
\begin{equation}
B^{H}(t)=c_{H}\int_{-\infty}^{t}\left((t-s)^{H-1/2}-(-s)_{+}^{H-1/2}\right)\, dB(s),\label{eq:c1}
\end{equation}
 where $\{B(t)\}_{t\in\R}$ is a standard Brownian motion. Then
 \begin{equation}
c_{H}=\sqrt{2H}\,2^{H}B\left(1-H,H+1/2\right)^{-1/2}\,,\label{eq:c2}
\end{equation}
 where $B(a,b) = \int_0^1 x^{a-1}(1-x)^{b-1}\, dx$ is the usual beta function.
 \end{lemma}

\textbf{Proof.} Since $\text{Var}(B^{H}(1))=1$ we get
\begin{equation}
c_{H}=\left\{ \int_{0}^{\infty}\left((1+x)^{H-1/2}-x^{H-1/2}\right)^{2}\, dx+\frac{1}{2H}\right\} ^{-1/2}.\label{eq:c3}
\end{equation}
 Put
 \[
I=\int_{0}^{\infty}\left((1+x)^{H-1/2}-x^{H-1/2}\right)^{2}\, dx\,.
\]
 Then
 \begin{align*}
I & =\lim_{\mu\to0^{+}}\,\int_{0}^{\infty}\left((1+x)^{H-1/2}-x^{H-1/2}\right)^{2}e^{-\mu x}\, dx\\
 & =\lim_{\mu\to0^{+}}\left\{ (e^{\mu}+1)\mu^{-2H}\Gamma(2H)-e^{\mu}\mu^{-2H}\gamma(2H,\mu)-2\int_{0}^{\infty}(1+x)^{H-1/2}x^{H-1/2}e^{-\mu x}\, dx\right\} \\
 & =-\frac{1}{2H}+\lim_{\mu\to0^{+}}\left\{ 2e^{\mu/2}\mu^{-2H}\Gamma(2H)-2\int_{0}^{\infty}(1+x)^{H-1/2}x^{H-1/2}e^{-\mu x}\, dx\right\} \\
 & =-\frac{1}{2H}+\lim_{\mu\to0^{+}}\left\{ 2e^{\mu/2}\mu^{-2H}\Gamma(2H)-\frac{2}{\sqrt{\pi}}e^{\mu/2}\Gamma\big(H+\frac{1}{2}\big)\,\mu^{-H}K_{-H}\big(\frac{\mu}{2}\big)\right\} ,
 \end{align*}
 where $\gamma(z,x)$ and $K_{\nu}(z)$ are the incomplete gamma function
and modified Bessel function of the second kind, respectively. The
third equality uses the facts that $e^{\mu}\mu^{-2H}\gamma(2H,\mu)=\frac{1}{2H}$
+ o(1), and that $(e^{\mu}+1)\mu^{-2H}=2e^{\mu/2}\mu^{-2H}+o(1)$
for $H<1$, as $\mu\to0$. The forth equality applies formula 3.3838
in \cite{GR}.

Using the duplication formula $\Gamma(2H)=\frac{2^{2H-1}}{\sqrt{\pi}}\Gamma(H)\Gamma\big(H+\frac{1}{2}\big)$ (see \cite[formula 8.3351]{GR}), we get
 \begin{align*}
I & =-\frac{1}{2H}+\frac{1}{\sqrt{\pi}}\,\Gamma\big(H+\frac{1}{2}\big)\,\lim_{\mu\to0^{+}}\left\{ \mu^{-2H}2^{2H}\Gamma(H)-2\mu^{-H}\, K_{H}\big(\frac{\mu}{2}\big)\right\} .\label{eq:c4}
\end{align*}
 Since
 \[
\mu^{-2H}2^{2H}\Gamma(H)=\int_{0}^{\infty}x^{H-1}e^{-\frac{\mu^{2}}{4}x}\, dx  \quad \text{and}  \quad K_{\nu}(z)=\frac{1}{2}\big(\frac{z}{2}\big)^{\nu}\int_{0}^{\infty}t^{-\nu-1}e^{-t-\frac{z^{2}}{4t}}\, dt
\]
(see \cite[formula 3.4326]{GR}), we obtain
\begin{align}
I & =-\frac{1}{2H}+\frac{1}{\sqrt{\pi}}\,\Gamma\big(H+\frac{1}{2}\big)\,\lim_{\mu\to 0^{+}} \int_{0}^{\infty}x^{H-1}e^{-\frac{\mu^{2}}{4}x}(1-e^{-\frac{1}{4x}})\, dx \nonumber \\
 & = -\frac{1}{2H}+\frac{1}{\sqrt{\pi}}\,\Gamma\big(H+\frac{1}{2}\big)\,\int_{0}^{\infty}x^{H-1}(1-e^{-\frac{1}{4x}})\, dx \nonumber \\
 & = -\frac{1}{2H}+\frac{\Gamma(1-H)\Gamma\big(H+\frac{1}{2}\big)}{\sqrt{\pi}\,4^{H}H}. \label{eq:c5}
\end{align}
Combining \eqref{eq:c5} with \eqref{eq:c3} and using the well-known formula $B(x,y)=\dfrac{\Gamma(x)\Gamma(y)}{\Gamma(x+y)}$
(see, e.g., \cite[formula 8.3841]{GR}), we get \eqref{eq:c2}.
\hfill \qed

%\noindent Below we give a plot of $c_{H}$ as a function of $H\in(0,1)$.
%\bigskip{}
%\begin{center}
%\includegraphics[width=3in]{Cfunction1}
%\par\end{center}
%\bigskip{}

%\subsection*{Determinant of a Gaussian covariance}
%
%\begin{lemma}\label{res}
%Let $x_{1},\dots,x_{n}$ be vectors in an inner product space $(E,\langle\cdot,\cdot\rangle)$  
%and let $C=\left[\langle x_{i},x_{j}\rangle\right]_{1\le i,j\le n}$.
%Then
%\[
%\mathrm{det}(C)=\|x_{1}\|^{2}\,\|x_{2}-\mathrm{proj}_{x_{1}}(x_{2})\|^{2}\cdots\|x_{n}-\mathrm{proj}_{x_{1},\dots,x_{n-1}}(x_{n})\|^{2}
%\]
% where $\mathrm{proj}_{x_{1},\dots,x_{i-1}}(x_{i})$ denotes the orthogonal
%projection of $x_{i}$ onto  $\mathrm{span}\{x_{1},\dots,x_{i-1}\}$.
%\end{lemma}
%
%\textbf{Proof.} Notice that vectors $y_{i}=x_{i}-\text{proj}_{x_{1},\dots,x_{i-1}}(x_{i})$,
%$i=1,\dots,n$, are orthogonal and
%\[
%x_{i}=a_{i1}y_{1}+\cdots+a_{in}y_{n}
%\]
% for some $a_{ij}\in\R$ with $a_{ii}=1$ and $a_{ij}=0$ for $j>i$.
%Since
%\[
%\langle x_{i},x_{j}\rangle=\sum_{k=1}^{n}a_{ik}a_{jk}\|y_{k}\|^{2},
%\]
% we have $C=AA^{\top}$, where $A=\left[a_{ij}\|y_{j}\|\right]_{1\le i,j\le n}$
%is a lower triangular matrix with $\|y_{i}\|$'s on the diagonal.
%Hence
%\[
%\text{det}(C)=\text{det}(A)^{2}=\prod_{i=1}^{n}\|y_{i}\|^{2}.
%\]
%\hfill \qed

\bigskip	

\section*{Acknowledgments}

{\small The authors are very grateful to the referees for their careful reading of the manuscript and useful comments.}

\bibliographystyle{amsplain}

\end{document}